\documentclass[12pt,a4paper,reqno]{amsart} 
\usepackage{amssymb}
\usepackage{latexsym}
\usepackage{amsmath}
\usepackage{mathrsfs}
\usepackage[X2,T1]{fontenc}
\usepackage[applemac]{inputenc}
\usepackage{cite}
\usepackage{calc}                   

\newcommand{\dir}{{\operatorname{dir}}}
\newcommand{\Char}{\operatorname{Char}}
\newcommand{\scal}[2]{\langle #1,#2\rangle}
\newcommand{\rr}[1]{\mathbf R^{#1}}

\newcommand{\nm}[2]{\Vert #1\Vert _{#2}}

\newcommand{\op}{\operatorname{Op}}

\newcommand{\sets}[2]{\{ \, #1\, ;\, #2\, \} }

\newcommand{\ep}{\varepsilon}
\newcommand{\fy}{\varphi}
\newcommand{\cdo}{\, \cdot \, }
\newcommand{\supp}{\operatorname{supp}}

\newcommand{\eabs}[1]{\langle #1\rangle}     

\newcommand{\vrum}{\vspace{0.1cm}}

\newcommand{\back}[1]{\backslash {#1}}
\newcommand{\Op}{\operatorname{Op}}
\newcommand{\SG}{\operatorname{SG}}

\newcommand{\FB}{\mathscr F\! \mathscr B}

\newcommand{\WFF}{\operatorname{WF}}

\newcommand{\WF}[1]{\operatorname{WF}^{#1}}

\newcommand{\WFMB}{\operatorname{WF}_{M(\omega,\mathscr B)}}
\newcommand{\WFM}{\operatorname{WF}_{M(\omega,\mathscr B)}'}

\newcommand{\WFpB}{\operatorname{WF}^{\psi}_{M(\omega,\mathscr B)}}

\newcommand{\WFFBog}{\operatorname{WF}_{\FB_0(\omega)}}
\newcommand{\WFB}{\operatorname{WF}_{\FB(\omega)}'}

\newcommand{\MBo}{M(\omega,\mathscr B)}

\newcommand{\norm}[1]{\langle#1\rangle}

\newcommand{\capup}{\cap \pmb \cup}
\newcommand{\cupap}{\cup \pmb \cap}

\def\cB{{\mathcal B}}
\def\cC{{\mathcal C}}

\def\cS{{\mathscr S}}

\def\ccD{{\mathcal D}}

\newcommand{\WFgB}{\operatorname{WF}_{\cB}}
\newcommand{\WFgpB}{\operatorname{WF}^{\psi}_{\cB}}
\newcommand{\WFgeB}{\operatorname{WF}^{e}_{\cB}}
\newcommand{\WFgpeB}{\operatorname{WF}^{\psi e}_{\cB}}

\numberwithin{equation}{section}          

\newtheorem{thm}{Theorem}
\numberwithin{thm}{section}
\newtheorem*{tom}{\rubrik}
\newcommand{\rubrik}{}
\newtheorem{prop}[thm]{Proposition}
\newtheorem{cor}[thm]{Corollary}
\newtheorem{lemma}[thm]{Lemma}

\theoremstyle{definition}

\newtheorem{defn}[thm]{Definition}
\newtheorem{example}[thm]{Example}

\theoremstyle{remark}
\newtheorem{rem}[thm]{Remark}

\author{Sandro Coriasco}

\address{Dipartimento di Matematica, Università di Torino, Italy}

\email{sandro.coriasco@unito.it}

\author{Karoline Johansson}

\address{Department of Computer science, Physics and Mathematics, Linn{\ae}us University, V{\"a}xj{\"o}, Sweden\footnote{Former address: Department of Mathematics and Systems Engineering,
V{\"a}xj{\"o} University, Sweden}}

\email{karoline.johansson@lnu.se}

\author{Joachim Toft}

\address{Department of Computer science, Physics and Mathematics, Linn{\ae}us University, V{\"a}xj{\"o}, Sweden${}^1$}

\email{joachim.toft@lnu.se}

\title{Global wave-front sets of Banach, Fr{\'e}chet and Modulation space types, and pseudo-differential operators}

\keywords{Wave-front, Fourier, Banach space, Modulation space, Micro-local, pseudo-differential}

\subjclass[2000]{35A18,35S30,42B05,35H10}

\begin{document}

\setlength{\textwidth}{125mm}
\setlength{\textheight}{185mm}

\begin{abstract}
We introduce global wave-front sets $\WFgB (f)$, $f\in
\cS^\prime(\rr d)$, with respect to suitable Banach or Fr\'echet
spaces $\cB$. An important special case is given by the modulation spaces
$\cB=\MBo$, where $\omega$ is an appropriate weight function and $\mathscr B$
is a translation invariant Banach function space. We show that the
standard properties for known notions of wave-front set extend to
$\WFgB (f)$. In particular, we prove that microlocality and
microellipticity hold for a class of globally defined
pseudo-differential operators $\op_t(a)$, acting continuously
on the involved spaces.
\end{abstract}

\maketitle

\section{Introduction}\label{sec0}

\par

An important question for a linear operator $T$ is wether $T$ possess "convenient" invertibility properties. For example, if $T=1-\Delta$, where $\Delta$ is the Laplace operator, and $f\in \mathscr S'(\rr d)$ fulfills
\begin{equation}\label{problemform}
Tf = g,
\end{equation}
for some given $g\in \mathscr S(\rr d)$, then it follows by Fourier's inversion formula that $f\in \mathscr S(\rr d)$. Similar facts are true when $1-\Delta$ is replaced by any partial differential operator whose symbol satisfies certain hypoelliptic condition. Furthermore, if in addition the coefficients are constant, then $f$ can easily be computed by a similar Fourier argument as for $1-\Delta$.

\par

In the case that some of the coefficients are not constant, the question of finding $f$ is more complicated, because the Fourier technique can not be applied in such simple way. In this situation one may approximate $f$ by discretizing the problem, using Gabor analysis in combination with pseudo-differential calculus. In this setting one may assume that $f$ and $g$ belong
to appropriate modulation spaces which are possible to discretize in the
context of Gabor theory. This implies that \eqref{problemform} can be formulated in such way that $f$ and $g$ are replaced by sequences in appropriate $l^p-l^q$ spaces, and $T$ is replaced by an infinite matrix with "good" control of the involved coefficients. (See e.{\,}g. papers in \cite{FeiStr} for details.)

\par

In the case that the symbol for $T$ does not fulfill such hypoelliptic condition, the situation is even
more complicated, and \eqref{problemform} might not be sufficient to conclude that $f$ should belong to $\mathscr S(\rr d)$ when $g\in \mathscr S(\rr d)$. In this context it is interesting to know \emph{where} and \emph{how} $f$ fails to belong to $\mathscr S$. A convenient way to consider such questions is to use global wave-front sets of smooth types, introduced in \cite{Me:94} and, with a different
approach, given in \cite{CoMa}. In these papers, such wave-front sets are used for performing global and
detailed regularity investigations within the
theory of partial differential equations, for
pseudo-differential operators as well as for Fourier integral operators.
For example, convenient mapping
properties for pseudo-differential operators and Fourier integral operators 
of $\SG$-type on global
wave-front sets with respect to smoothness are established in
\cite{CoMa}. These properties are used to
describe propagation of singularities for a class of hyperbolic Cauchy
problems, associated with linear operators having coefficients growing
(at most) polynomially with respect to the space
variable (cf.  \cite[Theorem 5.4]{CoMa}). See also e.{\,}g.
\cite{Co, coriasco, CoRo, CaRo, Me:94, Parenti} for related
results and investigations. 

\medspace

In this paper we introduce global wave-front sets for temperate distributions
with respect to certain Banach or Fr\'echet spaces $\cB$.
In particular, we focus on the case of modulation spaces $\cB=M(\omega,\mathscr B)$,
parameterized with the translation invariant Banach function space
$\mathscr B$ and the weight function $\omega$. We refer to
these families of global wave-front sets as wave-front sets of $\cB$-type or 
modulation space type, respectively.
We also establish mapping properties for global wave-front sets under the action of a
broad class of pseudo-differential operators of $\SG$-type
(cf. e.{\,}g. \cite{Co, ES, Sc, Me:94, Parenti}). In this context, our
approach links the ``local'' analysis carried out in \cite{CJT1,PTT1},
with the ``global'' analogue treated, e.{\,}g., in \cite{CoMa}. In
particular, we recover the microlocality and
microellipticity properties that hold for wave-front sets of Sobolev
type introduced by H{\"o}rmander \cite{Ho1}, and classical wave-front
sets with respect to smoothness (cf. Sections 8.1 and 8.2 in
\cite{Ho2}), as well as for wave-front sets of Banach function type
in \cite{CJT1}, and wave-front sets with respect to $\mathscr S$ and
$H^2_{s_1,s_2}$ in \cite{CoMa} and \cite{Me:94}. Here $H^2_{s_1,s_2}$ is denoted by
$H^{s_1,s_2}$ in \cite{CoMa, Me:94}.
We also extend the mapping
properties in \cite{CoMa} to a broader class of
pseudo-differential operators, which contains those considered in \cite{CoMa}.

\par

The family of global wave-front sets of $\cB$-type
contains all the wave-front sets mentioned here above. This is a
consequence of the fact that the family of modulation spaces (which is a subfamily of all admissibe $\cB$) is
locally the same as the family of Fourier BF-spaces in \cite{CJT1},
and contains the $H^2_{s_1,s_2}$ spaces in \cite{Me:94} and the Sobolev
spaces in \cite{Hrm-nonlin}. Furthermore, for an appropriate family of
modulation spaces, their union and intersection are equal to
$\mathscr S'$ and $\mathscr S$ respectively.

\par

Here we recall that modulation spaces were introduced in \cite{F1} by
Feichtinger and further developed and generalized in
\cite{Feichtinger3, Feichtinger4, Feichtinger5, Grochenig0a}. The
modulation space $\MBo$, where $\omega$ denotes a weight and $\mathscr
B$ denotes a
(translation) invariant Banach function space on phase (or
time-frequency shift) space $\rr {2d}$, appears as the set of temperate
(ultra-) distributions whose short-time Fourier transforms multiplied with $\omega$ belong to $\mathscr B$
(cf. \cite{F1,F2,Feichtinger3,Feichtinger4}).

\par

In order to be more specific, let $\mathcal B$ be an appropriate
Banach space, or, more generally an appropriate Fr\'echet space of
distributions, and let $f$ be a temperate distribution on $\rr d$. In analogy with \cite{CoMa}, the
global wave-front set $\WFF _{\mathcal B}(f)$ of $f$ with respect to $\mathcal B$, is the union of three
components
\begin{equation}\label{globalWFcomponents}
\WF \psi _{\mathcal B}(f),\quad \WF e _{\mathcal B} (f)\quad
\text{and}\quad \WF {\psi e} _{\mathcal B}(f),
\end{equation}
where $\WF \psi _{\mathcal B}(f)$ agrees with the ``local'' wave-front
set which is carefully explained in \cite{CJT1,PTT1}. In the case
$\mathcal B=\mathscr S$, then $\WF \psi
_{\mathcal B}(f)$ is the same
as the ``classical'' wave-front set (cf. e.{\,}g. \cite[Sections
8.1--8.3]{Ho1}).
For convenience we usually write $\WF \psi (f)$, $\WF e
(f)$ and $\WF {\psi e} (f)$ instead of $\WF \psi _{\mathscr S}(f)$,
$\WF e _{\mathscr S} (f)$ and $\WF {\psi e} _{\mathscr S}(f)$,
respectively.
Note here that $\mathscr S$ is locally the same as
$C^\infty$, which implies that $\WF \psi _{\mathcal B}(f)$ is invariant under
the choices $\mathcal B=\mathscr S$ and $\mathcal B=C^\infty$.
We refer to the wave-front sets in
\eqref{globalWFcomponents} as the wave-front sets (with respect to
$\mathcal B$) for $f$ of $\psi$-type, $e$-type and $\psi e$-type,
respectively.

\par

Roughly speaking, $\WF \psi _{\mathcal B} (f)$ informs
about local singularities with respect to $\mathcal B$ and the
directions for their propagation. The set $\WF e _{\mathcal B}
(f)$ is essentially the same as $\WF \psi _{\mathcal B_0}(\widehat f)$
and informs about the directions were the size of $f$ fails to belong to
$\mathcal B$ near infinity. Here $\widehat f$ is the Fourier transform for $f$,
and $\mathcal B_0$ is a Banach or Frech\'et space related to $\mathcal
B$. Finally $\WF {\psi e} _{\mathcal B}(f)$ informs about those
directions were $f$ oscillates heavily at infinity compared to its size.

\par

Therefore, it might not be surprising that, if $\mathcal B$ is
appropriate, then the union $\WF e_{\mathcal
B}(f)\bigcup \WF {\psi e}_{\mathcal B} (f)$, the so-called "exit
component" (cf. \cite{CoMa}), explains where $f$, far away from origin,
fails to belong to $\mathcal B$. Taking into account
that $\WF \psi _{\mathcal B}(f)=\emptyset$, if and only if $f$ locally
belongs to $\mathcal B$, it follows that the global wave-front set
$\WFF _{\mathcal B}(f)$ fulfills
\begin{equation}\label{WFequiv}
\WFF _{\mathcal B}(f) =\emptyset \quad \Longleftrightarrow \quad f\in
{\mathcal B}
\end{equation}
for such $\mathcal B$. This relation was established in
\cite{CoMa} when $\mathcal B=\mathscr S$ (cf. \cite[Theorem
3.8]{CoMa}). In Theorem \ref{mainthm4} we prove \eqref{WFequiv} for general admissible $\mathcal B$.

\medspace

Let $\op _t (a)$ be a pseudo-differential operator with appropriate symbol
$a$, and let $\Char ^\psi (a)$ denote the set of characteristic
points for $a$ (or $\op _t (a)$). A fundamental mapping property for wave-front sets is given by the embeddings
\begin{equation}\label{standardWFemb}
\WF \psi _{\mathcal C}(\op _t(a)f) \subseteq \WF \psi _{\mathcal B}(f)
\subseteq \WF \psi _{\mathcal C}(\op _t(a)f) \bigcup \Char ^\psi (a)
\end{equation}
when $\mathcal B=\mathcal C =\mathscr S$ (cf. \cite[Section
18.1]{Ho1}).
The latter inclusions were extended to involve wave-front
sets of modulation space type in \cite{CJT1,PTT1}.

\par

Furthermore, if
$\mathcal B=\mathcal C =\mathscr S$, then, for
pseudo-differential operators of $\SG$-type, the inclusions in
\eqref{standardWFemb} were extended in \cite{CoMa} to include embeddings
involving wave-front sets of $e$- and $\psi e$-types, with modified sets
of characteristic points (see \cite[Proposition 3.6]{CoMa}).

\par

We extend and improve the latter inclusions in different ways. For example,  if we let $\eabs x=(1+|x|^2)^{1/2}$ when $x\in \rr d$, as usual, then we have the following:
\begin{itemize}
\item Our result holds for a broader class of pseudo-differential operators of $\SG$-type. More precisely, the symbol $a$ here above should belong to the set denoted by $\SG^{(\omega)}_{r,\rho}(\rr {2d})$ or $\SG^{(\omega)}_{r,\rho}$, for $r,\rho\ge 0$ and a weight function $\omega$. This set consists of all smooth functions $a$ on $\rr {2d}$ such that the global estimates
\begin{equation}\label{Somegadef}
	|D_x^\alpha D_\xi ^\beta a(x,\xi)|\leq C_{\alpha,\beta}\, \omega(x,\xi)\eabs
        x^{-r|\alpha|}\eabs \xi ^{-\rho|\beta|}, \qquad
        x,\xi\in \rr d,
\end{equation}
hold for suitable constants $C_{\alpha,\beta} \geq 0$, which only depend on the multi-indices $\alpha$ and $\beta$.

\vrum

\item Our sets of characteristic points are smaller compared to those in \cite{CoMa}.

\vrum

\item Our results hold for wave-front sets of $\mathcal B$- and $\mathcal
C$-types, when $\mathcal B$ and $\mathcal C$ are related to each other in an appropriate way. This family of wave-front sets contain the wave-front sets of Schwartz-type, and weighted Sobolev type, or more generally, modulation space type, as special cases.
\end{itemize}

\par
 
Moreover, we extend the theory and introduce global wave-front sets with respect to sequences $(\cB _j)$ and $(\cB _{j,k})$ of admissible spaces $\cB _j$ and $\cB _{j,k}$, and prove that all results carry over to this situation. For example, \eqref{WFequiv} and \eqref{standardWFemb} hold in the more general situation. By choosing $\cB_j$ and $\cB_{j,k}$ in appropriate ways, we show how one can define wave-front sets with respect to certain vector spaces $Q$ and $Q_0$ of smooth functions, which together with their derivatives should be bounded by certain polynomials. (Cf. Remark \ref{remhormwfsets2} for strict definition of $Q$ and $Q_0$.) Furthermore it follows from our investigations that \eqref{WFequiv} and \eqref{standardWFemb} hold when $\cB =Q$ or $\cB =Q_0$. 

\vrum

The results in the paper
expand the theory of modulation spaces in relation
to their interplay with pseudo-differential operators.
These topics have undergone an intense development
in the recent years, involving, e.{\,}g., the study of PDEs, 
as well as the boundedness properties
of pseudo-differential and Fourier integral operators,
on these
function and distribution spaces (see, e.g., \cite{BaHu, CoNi, CNR, CoZu}
and the references quoted therein).
Here we lay the foundations of the theory concerning the
spectral analysis of singularities of solutions of
PDEs in this context, in global
terms. That is, we take
care of the local properties of the solutions in
\cite{CJT1}, as well as, at the same time, of their
behaviour at infinity, in the spirit of \cite{CoMa}.
Especially, we here consider hypoelliptic and elliptic
operators.

\par

We remark that the notion of
wave-front set studied below generalizes
(and includes) other similar concepts, previously
examined in separate environments. Here we treat
them in a unified manner.
Recently there have been several interesting results on local and global
propagation of singularities for Schr{\"o}dinger equations in the context
of analytic function theory 
(see, e.g., \cite{MaNaSo:09,Na:05,RZ,Sj}). These investigations are partially
based on properties on (local) analytic wave-front sets. Here we emphasize
that we always have $\WFF _{\mathcal{B}}(f)\subseteq \WFF _{\mathscr{S}}(f)$, 
while opposite inclusion is valid when $\WFF _{\mathcal{B}}$ is replaced by
the analytic wave-front set. In particular, we are here considering other types
of singularities compared to  \cite{MaNaSo:09,Na:05,RZ,Sj}.

\par

We plan to further expand these studies
in forthcoming papers, with the aim of obtaining
analogues of known theorems which describe
wave-front sets in terms of unions of bicharacteristics, and
of analyzing the propagation of singularities for
evolution PDEs of strictly hyperbolic type and of Schr{\"o}dinger type.

\medspace

The paper is organized as follows. In Section \ref{sec1} we recall the
definition and basic properties of pseudo-differential operators,
translation invariant Banach function spaces and modulation
spaces. Here we also define three types of sets of characteristic
points and show some properties for them. One of these characteristic sets
coincide with the one defined in \cite{CJT1}, but all of them differ
from the corresponding ones considered in \cite{CaRo} and \cite{CoMa}.

\par

In Section \ref{sec2} we define the wave-front sets of $\cB$-type $\WFgB(f)$.
Furthermore, we prove that $\WFMB^\psi(f)$
and $\WFFBog^\psi(f)$ coincide with the wave-front sets defined in
\cite{CJT1}, when $M(\omega ,\mathscr B)$ is locally the same as the
Fourier BF-space $\FB _0(\omega )$. The  remaining part of
the section is devoted to the proof of a
relation between the wave-front sets of $\cB$-type and the
sets of characteristic points. Here we also prove \eqref{WFequiv} when
$\mathcal B= \MBo$.

\par

The Sections \ref{sec3} and \ref{sec4} are devoted to mapping properties
for pseudo-differential operators in the context of these wave-front
sets. The properties are given and proved in Section \ref{sec3}, and
similar properties are stated in Section \ref{sec4} for wave-front sets
with respect to sequences of modulation spaces. Examples of 
applications of our results are given in Section \ref{sec:ex} and at the end
of Section \ref{sec4}. Finally, in Appendix A
we prove some of the properties for the
sets of characteristic points stated in Section \ref{sec1}.

\par

\subsection*{Acknowledgement} 

\par

We are grateful to Professors Feichtinger and Rodino, for valuable advice and constructive critisism.

\par

\section{Preliminaries}\label{sec1}

\par

In this section we recall some notation and basic results. The proofs
are in general omitted. In what
follows we let $\Gamma$ denote an open cone in $\rr d\setminus 0$. (The vertices
of the cones are always origin.) If
$\xi \in \rr d\setminus 0$ is fixed, then an open cone which contains
$\xi $ is sometimes denoted by $\Gamma_\xi$.

\par

Let $\omega$ and $v$ be positive measurable functions
on $\rr d$. Then $\omega$ is called $v$-moderate if
\begin{equation}\label{moderate}
\omega (x+y) \leq C\omega (x)v(y)
\end{equation}
for some constant $C$ which is independent of $x,y\in \rr d$. If $v$
in \eqref{moderate} can be chosen as a polynomial, then $\omega$ is
called polynomially moderated. We let $\mathscr P(\rr d)$ be the set
of all polynomially moderated functions on $\rr d$. If $\omega (x,\xi
)\in \mathscr P(\rr {2d})$ is constant with respect to the
$x$-variable or the $\xi$-variable, then we sometimes write $\omega
(\xi )$, respectively $\omega (x)$, instead of $\omega (x,\xi )$. In
this case we consider
$\omega$ as an element in $\mathscr P(\rr {2d})$ or in $\mathscr P(\rr
d)$ depending on the situation. We say that $v$ is submultiplicative
if \eqref{moderate} holds for $\omega=v$. For convenience we assume
that all submultiplicative weights are even. In what follows we always let
$v$ and $v_j$ stand for submultiplicative weights, if nothing else is stated.

\par

We also need to consider classes of weight functions, related to
$\mathscr P$. More precisely, we let $\mathscr P_0(\rr d)$ be the set
of all $\omega \in \mathscr P(\rr d)\bigcap C^\infty (\rr d)$ such
that $\partial ^\alpha \omega /\omega \in L^\infty$ for all
multi-indices $\alpha$. By Lemma 1.2 in \cite {To8} it follows
that for each $\omega \in \mathscr P(\rr d)$, there is an element
$\omega _0\in \mathscr P_0(\rr d)$ which is equivalent to $\omega$ in
the sense that
\begin{equation}\label{onestar}
C^{-1}\omega _0\le \omega \le C\omega _0,
\end{equation}
for some positive constant $C$.

\par

We need some more conditions for the weights in $\SG ^{(\omega
)}_{r,\rho }(\rr {2d})$ (cf. \eqref{Somegadef}). More precisely let $r,\rho \ge
0$. Then we let $\mathscr
P_{r,\rho}(\rr {2d})$ be the set of
all $\omega (x,\xi )$ in $\mathscr P(\rr {2d})\bigcap C^\infty (\rr
{2d})$ such that
$$
\eabs x^{r|\alpha |}\eabs \xi ^{\rho |\beta |}\frac {\partial ^\alpha
_x\partial ^\beta _\xi \omega (x,\xi )}{\omega (x,\xi )}\in L^\infty
(\rr {2d}),
$$
for every multi-indices $\alpha$ and $\beta$. Note that $\mathscr
P_{r,\rho} $ is different here compared to \cite{CJT1}. Note also that
in contrast to
$\mathscr P_0$, we do not have an equivalence between $\mathscr
P_{r,\rho}$ and $\mathscr P$ when $r>0$ or $\rho >0$ in the sense of
\eqref{onestar}. On the other hand, if $s, t\in
\mathbf R$ and $r, \rho \in [0,1]$, then $\mathscr P_{r,\rho} (\rr {2d})$
contains $\omega (x,\xi )=\eabs x ^t\eabs \xi ^s$, which gives one of
the most important classes in the applications.

\medspace

The Fourier transform $\mathscr F$ is the linear and continuous
mapping on $\mathscr S'(\rr d)$ which takes the form
$$
(\mathscr Ff)(\xi )= \widehat f(\xi ) \equiv (2\pi )^{-d/2}\int _{\rr
{d}} f(x)e^{-i\scal  x\xi }\, dx
$$
when $f\in L^1(\rr d)$. We recall that $\mathscr F$ is a homeomorphism
on $\mathscr S'(\rr d)$ which restricts to a homeomorphism on $\mathscr
S(\rr d)$ and to a unitary operator on $L^2(\rr d)$.

\medspace

Next we define the Banach function spaces (BF-spaces) and present some
useful properties.

\par

\begin{defn}\label{BFspaces}
Let $\mathscr B$ be a
Banach space which is continuously embedded in $L^1_{\mathrm{loc}}(\rr d)$, and let $v \in \mathscr P(\rr {d})$
be submultiplicative.
Then $\mathscr B$ is called a \emph{(translation) invariant
BF-space on $\rr d$} (with respect to $v$), if there is a constant $C$
such that the following conditions are fulfilled:
\begin{enumerate}
\item $\mathscr S(\rr d)\subseteq \mathscr
B\subseteq \mathscr S'(\rr d)$ (continuous embeddings);

\vrum

\item if $x\in \rr d$ and $f\in \mathscr B$, then $f(\cdot -x)\in
\mathscr B$, and
\begin{equation}\label{translmultprop1}
\nm {f(\cdot-x)}{\mathscr B}\le Cv(x)\nm {f}{\mathscr B}\text ;
\end{equation}

\vrum

\item if  $f,g\in L^1_{\mathrm{loc}}(\rr d)$ satisfy $g\in \mathscr B$
and $|f| \le |g|$ almost everywhere, then $f\in \mathscr B$ and
$$
\nm f{\mathscr B}\le C\nm g{\mathscr B}\text .
$$
\end{enumerate}
\end{defn}

\par

Assume that $\mathscr B$ is a translation invariant BF-space. If $f\in
\mathscr B$ and $h\in L^\infty$, then it follows from (3) in
Definition \ref{BFspaces} that $f\cdot h\in \mathscr B$ and
\begin{equation}\label{multprop}
\nm {f\cdot h}{\mathscr B}\le C\nm f{\mathscr B}\nm h{L^\infty}.
\end{equation}

\par

Also let  $\omega \in \mathscr P(\rr {d})$. Then we let the Fourier
BF-space $\FB {(\omega )}$ be the set of all $f\in \mathscr S'(\rr d)$
such that $\xi \mapsto \widehat f(\xi )\omega (\xi )$ belongs to
$\mathscr B$. It follows that $\FB(\omega )$ is a Banach space under
the norm
\begin{equation}\label{FLnorm}
\nm f{\FB {(\omega )}}\equiv \nm {\widehat
f\, \omega }{\mathscr B}.
\end{equation}

\par

\begin{rem}\label{refBF-rum}
Definition \ref{BFspaces} is not the standard definition of Banach function spaces, and there are several approaches to the theory of such spaces. We refer to \cite{Fe1.5,LZ} and the references therein for a thorough review of the theory of such spaces.
\end{rem}

\par

\begin{rem}
In many situations it is convenient to permit an $x$ dependence for
the weight $\omega$ in the definition of Fourier BF-spaces.
More precisely, for each $\omega (x,\xi ) \in \mathscr P(\rr {2d})$
and each translation invariant
BF-space $\mathscr B$ on $\rr d$, we let $\FB(\omega)$ be the set of
all $f \in \mathscr S'(\rr d)$ such that
$$
\nm{f}{\FB(\omega)} = \nm{f}{\FB(\omega),x} \equiv \nm{\widehat f \, \omega(x,
\cdot )}{\mathscr B}
$$
is finite. Since $\omega$ is $v$-moderate for some $v \in \mathscr
P(\rr {2d})$ it follows that
different choices of $x$ give rise to equivalent norms. Therefore the
condition $\nm{f}{\FB(\omega)} < \infty$ is independent of $x$, and it
follows that $\FB(\omega)$ is independent of $x$ although
$\nm{\cdot}{\FB(\omega)}$ might depend on $x$.
\end{rem}

\par

If $\mathscr B$ is an invariant BF-space with respect to $v$, then we always assume that $\mathscr B$ fulfills the property that the convolution map $(f,g)\mapsto f*g$ is sequentially continuous from $\mathscr B\times C_0^\infty$ to $\mathscr B$. By Minkowski's inequality we get
\begin{multline*}
\nm {f*g}{\mathscr B} =\Big \Vert \int f(\cdo -y)g(y)\, dy \Big \Vert
_{\mathscr B}
\\[1ex]
\le \int \nm { f(\cdo-y)}{\mathscr B}|g(y)|\, dy \le C\int \nm {
f}{\mathscr B}|g(y)v(y)|\, dy = C\nm { f}{\mathscr B}\nm
g{L^1_{(v)}}.
\end{multline*}
when $f\in \mathscr B$ and $g\in  C_0^\infty$.
Since $C_0^\infty$ is dense in $L^1_{(v)}$ it follows that the convolution map
$(f,g)\mapsto f*g$ from $\mathscr B\times C_0^\infty$ to $\mathscr B$ extends
uniquely to a continuous mapping from $\mathscr B\times L^1_{(v)}$ to $\mathscr B$.

\medspace

In order to define modulation spaces, we need to consider local
Fourier transforms. First assume that $\phi \in \mathscr S(\rr
d)$. Then the \emph{short-time Fourier transform} of $f\in \mathscr
S(\rr d)$ with respect to (the window function) $\phi$ is defined by
\begin{equation}\label{stftformula}
V_\phi f(x,\xi ) = (2\pi )^{-d/2}\int _{\rr {d}} f(y)\overline {\phi
(y-x)}e^{-i\scal y\xi}\, dy.
\end{equation}
If instead $\phi \in \mathscr S'(\rr d)$ and  $f\in \mathscr
S'(\rr d)$, then
the short-time Fourier transform of $f$ with respect to $\phi$ is
defined by
\begin{equation}\tag*{(\ref{stftformula})$'$}
(V_\phi f) =\mathscr F_2F,\quad \text{where}\quad F(x,y)=(f\otimes
\overline \phi)(y,y-x).
\end{equation}
Here $\mathscr F_2F$ is the partial Fourier transform of
$F(x,y)\in \mathscr S'(\rr{2d})$ with respect to the $y$-variable.
The definition \eqref{stftformula}$'$ makes sense, since the mappings
$\mathscr F_2$ and $F\mapsto S^*F=F\circ S$ with $S(x,y)=(y,y-x)$ are homeomorphisms
on $\mathscr S'(\rr {2d})$. It is obvious that \eqref{stftformula} and \eqref{stftformula}$'$
agree when $f,\phi \in \mathscr S(\rr d)$.

\par

In the following lemma we recall some
general continuity properties of the short-time Fourier transform. We
omit the proof since the result can be found,  e.{\,}g., in \cite{Fo}.

\par

\begin{lemma}\label{stftlemma}
Let $T$ be the mapping from $\mathscr S'(\rr d)\times \mathscr S'(\rr
d)$ to $\mathscr S'(\rr {2d})$ which is given by $T(f,\phi )=V_\phi f$,
and assume that $f,\phi \in\mathscr S'(\rr d)\setminus 0$. Then the
following is true:
\begin{enumerate}
\item $T$ restricts to a continuous map from $\mathscr S(\rr d)\times
\mathscr S(\rr d)$ to $\mathscr S(\rr {2d})$. Furthermore, $V_\phi f\in
\mathscr S(\rr {2d})$ if and only if $f,\phi \in \mathscr S(\rr {d})$;

\vrum

\item $T$ restricts to a continuous map from $L^2(\rr d)\times L^2(\rr
d)$ to $L^2(\rr {2d})$. Furthermore, $V_\phi f\in L^2(\rr {2d})$ if and
only if $f,\phi \in L^2(\rr {d})$;

\vrum

\item $T$ is a continuous map from $\mathscr S'(\rr d)\times \mathscr
S'(\rr d)$ to $\mathscr S'(\rr {2d})$, and restricts to a continuous
map from $\mathscr S'(\rr d)\times \mathscr S(\rr d)$ and from
$\mathscr S(\rr d)\times \mathscr S'(\rr d)$ to $\mathscr S'(\rr
{2d})\bigcap C^\infty (\rr {2d})$.
\end{enumerate}
\end{lemma}

\par

For Lemma \ref{stftproperties} below, we recall some
formulae for the short-time Fourier transform
which are important in time-frequency analysis (see
e.{\,}g. \cite{Gro-book}). For this reason it is convenient to let
$\widehat *$ be the twisted convolution on $L^1(\rr {2d})$,
defined by the formula
$$
(F\, \widehat *\, G)(x,\xi )=(2\pi )^{-d/2}\iint F(x-y,\xi -\eta
)G(y,\eta )e^{-i\scal {x-y}\eta}\, dyd\eta .
$$
By straightforward computations it follows that $\widehat *$
restricts to a continuous multiplication on $\mathscr S(\rr
{2d})$. Furthermore, the map $(F,G)\mapsto F\, \widehat *\, G$ from
$\mathscr S(\rr {2d})\times \mathscr S(\rr {2d})$ to $\mathscr S(\rr
{2d})$ extends uniquely to continuous mappings from $\mathscr S'(\rr
{2d})\times \mathscr S(\rr {2d})$ and $\mathscr S(\rr {2d})\times
\mathscr S'(\rr {2d})$ to $\mathscr S'(\rr {2d})\bigcap C^{\infty}(\rr
{2d})$.

\par

The following lemma is a straightforward consequence of Fourier's
inversion formula. The proof is therefore omitted. Here and in what
follows we let $\widetilde f(x)=\overline {f(-x)}$.

\par

\begin{lemma}\label{stftproperties}
Assume that $f\in \mathscr S'(\rr d)$ and that $\phi ,\phi _j\in
\mathscr S(\rr d)$ for $j=1,2,3$. Then the following is true:
\begin{enumerate}
\item $(V_\phi f)(x,\xi )= e^{-i\scal x\xi }(V_{\widehat \phi}\widehat
f)(\xi ,-x)$;

\vrum

\item $(V_{\phi _1}f)\widehat *(V_{\phi _2}\phi _3) = (\phi _3,\phi
_1)_{L^2}\cdot V_{\phi _2}f$;

\vrum

\item $\displaystyle {(V_{\phi _1\overline{\phi _2}}f)(x,\xi )=(2\pi
)^{-d/2} \int _{\rr {d}}(V_{\phi _1}f)(x,\xi -\eta )\widehat \phi
_2(\eta )e^{-i\scal x\eta}\, d\eta}$;

\vrum

\item $\displaystyle {(V_{\phi _1*{\widetilde \phi _2}}f)(x,\xi )=
  \int _{\rr {d}} (V_{\phi _1}f)(x-y,\xi ) \phi _2(y )\, dy}$.
\end{enumerate}
\end{lemma}

\medspace

Now assume that $\mathscr{B}$ is a translation invariant BF-space on
$\rr {2d}$, with respect to $v\in \mathscr P(\rr {2d})$. Also let
$\phi \in \mathscr{S}(\rr{d})\back{0}$ and let $\omega \in
\mathscr{P}(\rr{2d})$ be such that $\omega$ is $v$-moderate. The
\emph{modulation space} $M(\omega ,\mathscr B)$ consists of all $f\in
\mathscr{S}'(\rr{d})$ such that $V_{\phi}f\cdot \omega \in
\mathscr{B}$. We note that $M(\omega ,\mathscr B)$ is a Banach space
with the norm
\begin{equation}\label{modnorm}
\|f\|_{M(\omega ,\mathscr B)} \equiv \|V_{\phi} f
\omega\|_{\mathscr{B}}
\end{equation}
(cf. \cite{Feichtinger3}).

\par

\begin{rem}
Assume that $p,q\in [1,\infty]$, and let $L^{p,q}_{1}(\rr{2d})$ and
$L^{p,q}_{2}(\rr{2d})$ be the sets of all $F\in  L^1_{\mathrm{loc}}
(\rr{2d})$ such that
\begin{equation*}
\|F\|_{L^{p,q}_1}  \equiv \Big ( \int \Big( \int |F(x,\xi)|^p\,
dx\Big )^{q/p}\,d\xi \Big )^{1/q}
<\infty
\end{equation*}
and
\begin{equation*}
\|F\|_{L^{p,q}_2} \equiv \Big ( \int \Big ( \int |F(x,\xi)|^p\,
d\xi \Big )^{q/p}\, dx\Big )^{1/q}<\infty .
\end{equation*}
Then $M(\omega ,L^{p,q}_1(\rr{2d}))$ is equal to the usual modulation
space $M^{p,q}_{(\omega)}(\rr{d})$, and $M(\omega ,L^{p,q}_2(\rr{2d}))$ is
equal to the space $W^{p,q}_{(\omega)}(\rr{d})$, related to
Wiener-amalgam spaces (cf. \cite{F1,Feichtinger6,Feichtinger3,Gro-book}).
\end{rem}

\par

For notational convenience we set $M^p_{(\omega )} = M^{p,p}_{(\omega
)}= W^{p,p}_{(\omega )}$. Furthermore, if $\omega =1$, then we write
$M^{p,q}$, $M^p$ and $W^{p,q}$ instead of $M^{p,q}_{(\omega )}$,
$M^p_{(\omega )}$ and $W^{p,q}_{(\omega )}$ respectively.

\par

In the following proposition we list some important properties for
modulation spaces. We refer to \cite{Gro-book} for the proof.

\par

\begin{prop}\label{modproperties}
Let $\omega _0, v_0,v\in\mathscr{P}(\rr{2d})$ be such that $v$ and
$v_0$ are submultiplicative, and  $\omega _0$ is $v_0$-moderate. Also
let $\mathscr B$ be a translation invariant BF-space on $\rr
{2d} $ with respect to $v$ and $f\in \mathscr S'(\rr d)$. Then the
following is true:
\begin{enumerate}
\item if $\phi\in M^1_{(v_0v)}(\rr{d})\back{0}$, then $f\in M(\omega
,\mathscr B)$ if and only if $V_\phi f\cdot \omega \in \mathscr
B$. Furthermore, \eqref{modnorm} defines a norm on $M(\omega
,\mathscr B)$, and different choices of $\phi$ give rise to equivalent
norms;

\vrum

\item $\mathscr S(\rr d)\subseteq M^1_{(v_0v)}(\rr d)\subseteq
M(\omega ,\mathscr{B})\subseteq M^{\infty}_{(1/(v_0v))}(\rr
d)\subseteq \mathscr S'(\rr d)$.
\end{enumerate}
\end{prop}

\par
Proposition \ref{modproperties}, (1) allows us to be rather vague about
the choice of $\phi \in  M^1_{(v)}\setminus 0$ in
\eqref{modnorm}. For example, if $C>0$ is a constant and $\Omega$ is a
subset of $\mathscr S'$, then $\nm a{M(\omega,\mathscr B )}\le C$ for
every $a\in \Omega$, means that the inequality holds for some choice
of $\phi \in  M^1_{(v)}\setminus 0$ and every $a\in
\Omega$. Evidently, for any other choice
of $\phi \in  M^1_{(v)}\setminus 0$, a similar inequality is true
although $C$ may have to be replaced by a larger constant, if
necessary.

\par

\begin{rem}
Several important spaces agree with certain modulation spaces. In fact, let $s,t\in \mathbf R$. Then the following is true:
\begin{enumerate}
\item if $\omega _0(x,\xi )=\eabs \xi ^s$ ($\omega _0(x,\xi )=\eabs x^t$), then $M^2_{(\omega _0)}(\rr d)$ is equal to the Sobolev space $H^2_s(\rr d)$ (the weighted $L^2$-space $L^2_t(\rr d)$);

\vrum

\item if $\omega _0(x,\xi )=\eabs x^t\eabs \xi ^s$, then $M^2_{(\omega _0)}(\rr d)$ is equal to the weighted
Sobolev space $H^2_{s,t}(\rr d)$ in \cite{CoMa,Me:94};
\end{enumerate}
\end{rem}

\par

We recall that Fourier BF-spaces and modulation spaces are locally the
same. In fact, let $\fy \in \mathscr S(\rr d)\back 0$, $\mathscr B$ be
a translation invariant BF-space on $\rr {2d}$, $\omega \in \mathscr
P(\rr {2d})$ and set $\omega _0 (\xi) =\omega (x _0,\xi)$ for some
fixed $x _0\in \rr d$. Then
\begin{equation}\label{B0def}
\mathscr{B}_0 \equiv \sets {f\in \mathscr
S'(\rr d)}{\fy \otimes f \in \mathscr B}
\end{equation}
is a translation invariant BF-space on $\rr d$ under the norm $\nm
f{\mathscr B_0}\equiv \nm {\fy \otimes f}{\mathscr B}$. The space
$\mathscr B_0$ is independent of $\fy \in \mathscr S(\rr d)\back 0$,
and different choices of $\fy$ gives rise of equivalent
norms. Furthermore
\begin{equation}\label{compsupps}
M(\omega ,\mathscr B)\bigcap \mathscr E'(\rr d) = \FB _0(\omega _0)\bigcap
\mathscr E'(\rr d)
\end{equation}
(cf. \cite{CJT1,RSTT}).

\medspace

Next we recall some facts in Chapter XVIII in \cite {Ho1}
concerning pseudo-differential operators. Let $a\in
\mathscr S(\rr {2d})$, and $t\in \mathbf R$ be fixed. Then
the pseudo-differential operator $\op _t(a)$ is the linear and
continuous operator on $\mathscr S(\rr d)$ defined by the formula
\begin{equation}\label{e0.5}
(\op _t(a)f)(x)
=
(2\pi ) ^{-d}\iint a((1-t)x+ty,\xi )f(y)e^{i\scal {x-y}\xi }\,
dyd\xi .
\end{equation}
For
general $a\in \mathscr S'(\rr {2d})$, the pseudo-differential
operator $\op _t(a)$ is defined as the continuous operator from
$\mathscr S(\rr d)$ to $\mathscr S'(\rr d)$ with distribution
kernel
\begin{equation}\label{weylkernel}
K_{t,a}(x,y)=(2\pi )^{-d/2}(\mathscr F_2^{-1}a)((1-t)x+ty,x-y).
\end{equation}
This definition makes sense, since $\mathscr F_2$ and the map
$$
F\mapsto F\circ S_t \;\text{ with }\; S_t(x,y)=((1-t)x+ty,x-y),
$$
are homeomorphisms on $\mathscr
S'(\rr {2d})$. We also note that the latter definition of $\op _t(a)$
agrees with the operator in \eqref{e0.5} when $a\in \mathscr S(\rr
{2d})$.

\par

If  $t=0$, then $\op _t(a)$ is the Kohn-Nirenberg
representation $\op (a)=a(x,D)$, and if $t=1/2$, then $\op _t(a)$ is the Weyl quantization. 

\par

Let $a\in \mathscr S'(\rr {2d})$ and $s,t\in
\mathbf R$. Then there is a unique $b\in \mathscr S'(\rr {2d})$ such
that $\op _s(a)=\op _t(b)$. By straightforward applications of
Fourier's inversion  formula, it follows that
\begin{equation}\label{pseudorelation}
\op _s(a)= \op _t(b) \quad \Longleftrightarrow \quad b(x,\xi
)=e^{i(t-s)\scal
{D_x}{D_\xi}}a(x,\xi)
\end{equation}
(cf. Section 18.5 in \cite{Ho1}).

\par

Next we discuss our symbol classes. Let $m,\mu ,r, \rho \in \mathbf R$
be fixed. Then $\SG
^{m,\mu }_{r,\rho}(\rr {2d})$ is the set of all $a\in
C^\infty (\rr {2d})$ such that for each pair of multi-indices
$\alpha$ and $\beta$, there is a constant $C_{\alpha ,\beta}$ such
that
$$
|D _x^\alpha D _\xi ^\beta a(x,\xi )|\le
C_{\alpha ,\beta }\eabs x^{m-r|\alpha|}\eabs \xi ^{\mu -\rho |\beta
|}.
$$
Usually we assume that $r,\rho \ge 0$ and that $\rho +r >0$.

\par

More generally, assume that $\omega \in \mathscr P_{r ,\rho}
(\rr {2d})$. Then we recall from the introduction that
$\SG _{r,\rho}^{(\omega )}(\rr {2d})$ consists of all $a\in
C^\infty (\rr {2d})$ such that for each pair of multi-indices
$\alpha$ and $\beta$, there are constants $C_{\alpha ,\beta}$ such
that \eqref{Somegadef} holds. We note that
\begin{equation}\label{SG}
\SG _{r,\rho}^{(\omega )}(\rr {2d})=S(\omega ,g_{r,\rho}
),
\end{equation}
when $g=g_{r,\rho}$ is the Riemannian metric on $\rr {2d}$,
defined by the formula
\begin{equation}\label{riemannianmetric}
\big (g_{r,\rho}\big )_{(y,\eta )}(x,\xi ) =\eabs y ^{-2r}|x|^2 +\eabs
\eta ^{-2\rho}|\xi |^2
\end{equation}
(cf. Section 18.4--18.6 in \cite{Ho1}). Furthermore, $\SG ^{(\omega
)}_{r,\rho} =\SG ^{m,\mu }_{r,\rho}$ when $\omega (x,\xi )=\eabs
x^m\eabs \xi ^\mu $.

\par

The following result shows that pseudo-differential operators with
symbols in $\SG ^{(\omega )}_{r,\rho}$ behave well. Here and in what follows we set
\begin{equation}\label{sigmammy}
\sigma _m(x)=\eabs x^m\quad \text{and}\quad \sigma _{m,\mu}(x,\xi )=\eabs {x}^m\eabs \xi ^\mu 
\end{equation}
when $m,\mu \in \mathbf R$, as in Proposition \ref{psicontprop}.

\par

\begin{prop}\label{psicontprop}
Let $\mathscr B$ be a translation invariant BF-space on $\rr {2d}$,
$s,t\in \mathbf R$, $r,\rho \ge 0$, $\omega \in \mathscr P(\rr {2d})$,
$\omega _0\in \mathscr P_{r,\rho} (\rr {2d})$ and that $a,b\in
\mathscr S'(\rr {2d})$ are such that $\op _s(a)=\op _t(b)$. Then the
following is true:
\begin{enumerate}
\item $a\in \SG ^{(\omega _0)}_{r,\rho} (\rr {2d})$ if and only if
$b\in \SG ^{(\omega _0)}_{r,\rho} (\rr {2d})$, and then
$$
a-b\in \SG ^{(\omega _0/\sigma _{r,\rho})}_{r,\rho} (\rr {2d})\text ;
$$

\vrum

\item if $a\in \SG ^{(\omega _0)}_{r,\rho} (\rr {2d})$, then $\op
_t(a)$ is continuous on $\mathscr S(\rr d)$ and extends uniquely to a
continuous operator on $\mathscr S'(\rr d)$;

\vrum

\item if $a\in \SG ^{(\omega _0)}_{r,\rho} (\rr {2d})$, then $\op _t(a)$
is continuous from $M(\omega ,\mathscr B)$ to $M(\omega /\omega
_0,\mathscr B)$;

\vrum

\item there exist $a\in \SG ^{(\omega _0)}_{r,\rho} (\rr {2d})$ and $b\in \SG ^{(1/\omega _0)}_{r,\rho} (\rr {2d})$ such that for every choice of $\omega \in \mathscr P(\rr {2d})$ and every translation invariant BF-space $\mathscr B$ on $\rr {2d}$, the mappings
\begin{gather*}
\op _t(a)\, : \, \mathscr S(\rr d)\to \mathscr S(\rr d),\quad \op _t(a)\, : \, \mathscr S'(\rr d)\to \mathscr S'(\rr d)
\\[1ex]
\text{and}\quad \op _t(a)\, : \, M(\omega ,\mathscr B)\to M(\omega /\omega _0,\mathscr B).
\end{gather*}
are continuous bijections with inverses $\op _t(b)$.
\end{enumerate}
\end{prop}

\par

\begin{proof}
From the assumptions it follows that $g_{r,\rho}$ in
\eqref{riemannianmetric} is slowly varying, $\sigma$-temperate and
satisfies $g_{r,\rho}\le g_{r,\rho}^\sigma$, and that $\omega$ is
$g_{r,\rho}$-continuous and ($\sigma$,$g_{r,\rho}$)-temperate
(see Sections 18.4--18.6 in \cite{Ho1} for definitions).
The assertions (1) and (2) are now consequences of
Proposition 18.5.10 and Theorem 18.6.2 in \cite{Ho1}, and \eqref{SG}.

\par

Finally, (3) and (4) follow immediately from \cite[Theorem 3.2]{Toft36} and \cite[Theorem 2.1]{GT}. The
proof is complete.
\end{proof}

\par

We remark that explicit bijections of the form in Proposition \ref{psicontprop} (4) can be found in \cite[Section 3]{GT}.

\par

The following definition is motivated by Proposition \ref{psicontprop} (3) and (4).

\par

\begin{defn}\label{admspacesdef}
Let $r,\rho \in [0,1]$, $t\in \mathbf R$, $\cB$ be a Banach or Fr{\'e}chet space of distributions on $\rr d$ such that
$$
\mathscr S(\rr d)\subseteq \cB  \subseteq \mathscr S'(\rr d)
$$
with continuous embeddings. Then $\cB$ is called \emph{$\SG$-admissible (with respect to $r$, $\rho$ and $d$)} when $\op _t(a)$ maps $\cB$ continuously into itself, for every $a\in \SG ^{0,0}_{r,\rho}$.

\par

If $\cB$ and $\cC$ are $\SG$-admissible with respect to $r$, $\rho$ and $d$, and $\omega _0\in \mathscr P_{r,\rho}(\rr {2d})$, then the pair $(\cB ,\cC )$ is called \emph{$\SG$-ordered (with respect to $\omega _0$)}, when the mappings
$$
\op _t(a)\, :\, \cB \to \cC \quad \text{and}\quad \op _t(b)\, :\, \cC \to \cB
$$
are continuous for every $a\in \SG ^{(\omega _0)}_{r,\rho}(\rr {2d})$ and $b\in \SG ^{(1/\omega _0)}_{r,\rho}(\rr {2d})$.
\end{defn}

\par

\begin{example}
Let $r,\rho \in [0,1]$, $t\in \mathbf R$, $\omega \in \mathscr P(\rr {2d})$ and let $\mathscr B$ be a translation invariant BF-spaceon $\rr {2d}$. Then Proposition \ref{psicontprop} shows that $\mathscr S(\rr d)$, $\mathscr S'(\rr d)$ and $M(\omega ,\mathscr B)$ are $\SG$-admissible. Furthermore, if $\omega _0\in \mathscr P_{r,\rho }(\rr {2d})$, then the pairs
$$
(\mathscr S(\rr d) ,\mathscr S(\rr d) ),\quad (\mathscr S'(\rr d), \mathscr S'(\rr d))\quad \text {and}\quad (M(\omega ,\mathscr B),M(\omega /\omega _0,\mathscr B) )
$$
are $\SG$-ordered with respect to $\omega _0$.
\end{example}

\par

\begin{rem}\label{remSGadm}
Let $\omega _0\in \mathscr P_{r,\rho}(\rr {2d})$ and let $\cB$ be $\SG$-admissible with respect to $r$, $\rho$ and $d$. Then there is a unique $\SG$-admissible $\cC$ such that $(\cB ,\cC )$ is a $\SG$-ordered pair with respect to $\omega _0$.

\par

In fact, let $a$ be as in Proposition \ref{psicontprop} (4). Then $\cC$ is the image of $\cB$ under $\op _t(a)$. The details are left for the reader.
\end{rem}

\par

If $a\in \SG^{(\omega _0 )}_{r,\rho}(\rr {2d})$, then it follows from
the definitions that there is a constant $C>0$ such that
$$
|a(x,\xi )|\le C\omega _0 (x,\xi ).
$$
On the other hand, a necessary and sufficient condition for $a$ to be
invertible, in the sense that $1/a$ should be a symbol in
$\SG^{(1/\omega _0 )}_{r,\rho}(\rr {2d})$, is that for some constant
$c>0$ we have
\begin{equation}\label{invcond}
c \, \omega _0 (x,\xi )\le |a(x,\xi )|.
\end{equation}
A slightly relaxed condition is that \eqref{invcond} holds for some
constant $c>0$ and all points $(x,\xi )$, outside a compact set
$K\subseteq \rr {2d}$. In this case we say that $a$ is \emph{elliptic}
(with respect to $\omega _0$).

\par

In the following we discuss more local invertibility conditions for
symbols in $\SG ^{(\omega _0)}_{r,\rho} (\rr {2d})$ in terms of
sets of characteristic points of the involved symbols. We remark that
our definition of such sets is slightly
different comparing to \cite[Definition 18.1.5]{Ho1} and \cite{CoMa} in view of Remark
\ref{compchar} below.

\par

\begin{defn}\label{defchar}
Assume that  $r,\rho \ge 0$, $\omega _0\in \mathscr P _{r,\rho} (\rr
{2d})$ and that $a\in \SG ^{(\omega _0)}_{r,\rho}(\rr {2d})$.

\medspace

\begin{enumerate}

\item $a$ is called
\emph{$\psi$-invertible} with respect to $\omega _0$ at the
point $(x_0,\xi_0)\in \rr d\times (\rr d\back 0)$, if
there exist a neighbourhood $X$ of $x_0$, an open conical
neighbourhood $\Gamma$ of $\xi _0$  and positive constants $R$ and $c$
such that \eqref{invcond} holds for $x\in X$, $\xi\in \Gamma$ and
$|\xi|\ge R$.

\par
The point $(x_0,\xi_0)$ is called \emph{$\psi$-characteristic} for $a$
with respect to $\omega _0$ if $a$ is \emph{not} $\psi$-invertible
with respect to $\omega _0$ at $(x_0,\xi_0)$;

\vrum

\item  $a$ is called
\emph{$e$-invertible} with respect to $\omega _0$ at the
point $(x_0,\xi_0)\in (\rr d\setminus 0)\times \rr d$, if
there exist an open conical neighbourhood $\Gamma$  of $x_0$, a
neighbourhood $X$ of $\xi _0$ and positive constants $R$ and $c$ such
that \eqref{invcond} holds for $x\in \Gamma$, $|x|\ge R$ and  $\xi\in
X$.

\par

The point $(x_0,\xi_0)$ is called \emph{$e$-characteristic} for $a$ with
respect to $\omega _0$ if $a$ is \emph{not} $e$-invertible
with respect to $\omega _0$ at $(x_0,\xi_0)$;

\vrum

\item  $a$ is called
\emph{$\psi e$-invertible} with respect to $\omega _0$ at the
point $(x_0,\xi_0)\in (\rr d\setminus 0)\times (\rr d\setminus 0)$, if
there exist open conical neighbourhoods $\Gamma _1$  of $x_0$ and
$\Gamma _2$ of $\xi _0$, and positive constants $R$ and $c$
such that \eqref{invcond} holds for $x\in \Gamma _1$, $|x|\ge R$,
$\xi \in \Gamma _2$ and $|\xi |\ge R$.

\par
The point $(x_0,\xi_0)$ is called \emph{$\psi e$-characteristic} for
$a$ with respect to $\omega _0$ if $a$ is \emph{not} $\psi
e$-invertible with respect to $\omega _0$ at $(x_0,\xi_0)$.

\end{enumerate}
The set of characteristic points (the characteristic set), for a
symbol $a\in \SG^{(\omega _0)}_{r,\rho}(\rr {2d})$ with respect to $
\omega _0$, is denoted by
$$
\Char (a)=\Char _{(\omega _0)}(a)=\Char ^\psi_{(\omega
_0)}(a)\bigcup\Char^e _{(\omega _0)}(a)\bigcup\Char ^{\psi e}_{(\omega
_0)}(a),
$$
where the three components are the sets of points satisfying (1), (2)
and (3), respectively.
\end{defn}

\par

\begin{rem}
In the case $\omega _0=1$ we exclude the phrase "with respect to $\omega _0$" in
Definition \ref{defchar}. For example, $a\in
\SG^{0,0}_{r,\rho}(\rr{2d})$ is \emph{$\psi$-invertible} at
$(x_0,\xi_0) \in \rr d\times (\rr d\back 0)$ if $(x_0,\xi _0)\notin
\Char ^{\psi} _{(\omega _0)}(a)$ with $\omega _0=1$. This means that there
exist a neighbourhood $X$ of $x_0$, an open conical
neighbourhood $\Gamma$ of $\xi_0$ and $R,c>0$ such that
\eqref{invcond} holds for $\omega _0=1$, $x\in X$ and $\xi\in \Gamma$
satisfies $|\xi|\geq R$.
\end{rem}

\par

In the next definition we introduce different classes of cutoff
functions (see also Definition 1.9 in \cite{CJT1}).

\par

\begin{defn}\label{cuttdef}
Assume that $X\subseteq \rr d$ is open, $\Gamma \subseteq \rr
d\setminus 0$ is an open cone, $x_0\in X$ and that $\xi _0\in \Gamma
$.

\begin{enumerate}
\item A smooth function $\fy$ on $\rr d$ is called a \emph{cutoff
function} with respect to $x_0$ and $X$, if $0\le \fy \le 1$, $\fy \in
C_0^\infty (X)$ and $\fy =1$ in an open neighbourhood of $x_0$. The
set of cutoff functions with respect to $x_0$ and $X$ is denoted by
$\mathscr C_{x_0}(X)$.

\vrum

\item  A smooth function $\psi$ on $\rr d$ is called a
\emph{directional cutoff function} with respect to $\xi_0$ and
$\Gamma$, if there is a constant $R>0$ and open conical neighbourhood
$\Gamma _1\subseteq \Gamma$ of $\xi _0$ such that the following is
true:
\begin{itemize}
\item $0\le \psi \le 1$ and $\supp \psi \subseteq \Gamma$;

\vrum

\item  $\psi (t\xi )=\psi (\xi )$ when $t\ge 1$ and $|\xi |\ge R$;

\vrum

\item $\psi (\xi )=1$ when $\xi \in \Gamma _1$ and $|\xi |\ge R$.
\end{itemize}

\par

The set of directional cutoff functions with respect to $\xi _0$ and
$\Gamma$ is denoted by $\mathscr C^\dir  _{\xi _0}(\Gamma )$.
\end{enumerate}
\end{defn}

\par

\begin{rem}
For notational convenience, the open neighbourhood $X$ of $x_0$ and
the open conical neighbourhood $\Gamma$ of $\xi_0$ appearing in the
previous definition will sometimes be omitted in the sequel, and we
will simply write $\mathscr C_{x_0}$ and $\mathscr C^\dir  _{\xi _0}$,
respectively.
\end{rem}

\par

\begin{rem}\label{psiinvremark}
Let $x_1,\xi_1\in \rr d$, $x_2,\xi_2 \in \rr d\back 0$, $\fy _1\in
\mathscr C _{x _1}(\rr d)$, $\fy _2\in \mathscr C _{\xi _1}(\rr d)$,
$\psi _1\in \mathscr C ^{\dir} _{x _2}(\rr d \back 0)$ and $\psi _2\in
\mathscr C ^{\dir} _{\xi _2}(\rr d \back 0)$. Then
\begin{equation}
c_1 = \fy _1\otimes \psi _2, \quad  c_2 = \psi _1\otimes \fy _2 \quad
\text{and} \quad c_3 = \psi _1\otimes \psi _2
\end{equation}
belong to $\SG^{0,0}_{1,1}(\rr {2d})$ and are $\psi$-invertible,
$e$-invertible and $\psi e$-invertible, respectively.
\end{rem}

\par

In the next three propositions we show that $\op _t(a)$ for $t\in \mathbf R$ satisfies
convenient invertibility properties of the form
\begin{equation}\label{locinvop}
\op _t(a)\op _t(b) = \op _t(c) + \op _t(h),
\end{equation}
outside the set of characteristic
points for a symbol $a$. Here $\op _t(b)$, $\op _t(c)$ and $\op _t(h)$ have
the roles of ''local inverse``, ''local identity`` and smoothing operators
respectively. From these propositions it also follows that our set of
characteristic points in Definition \ref{defchar} are related to those
in \cite{CoMa,Ho1}. 

\par

\begin{prop}\label{psicharequiv}
Let $r,\rho \in [0,1]$ be such that $\rho >0$, $\omega _0\in
\mathscr P_{r,\rho}(\rr {2d})$, $a\in \SG ^{(\omega _0)}_{r,\rho }(\rr
{2d})$, and let $(x_0,\xi _0)\in \rr d \times (\rr d \setminus
0)$. Then the following conditions are equivalent:

\medspace

\begin{enumerate}

\item $(x_0,\xi _0)\notin \Char _{(\omega _0)}^{\psi} (a)$;

\vrum

\item there is an element $c\in \SG^{0,0}_{r,\rho}$ which is
$\psi$-invertible at $(x_0,\xi _0)$, and an element $b\in
\SG^{(1/\omega _0)}_{r,\rho}$ such that $ab=c$;

\vrum

\item there is an element $c\in \SG^{0,0}_{r,\rho}$ which is
$\psi$-invertible at $(x_0,\xi _0)$, and elements  $h\in
\SG^{0,-\rho}_{r,\rho}$ and $b\in \SG^{(1/\omega _0)}_{r,\rho}$ such
that \eqref{locinvop} holds;

\vrum 

\item for each neighbourhood $X$ of $x_0$ and conical neighbourhood
$\Gamma $ of $\xi _0$, there is an element $c=\fy \otimes \psi$ where
$\fy \in \mathscr C_{x_0}(X)$ and $\psi \in \mathscr C_{\xi _0}^\dir
(\Gamma )$, and elements $h\in \mathscr S$ and $b\in \SG^{(1/\omega
_0)}_{r,\rho}$ such that \eqref{locinvop} holds. Furthermore, the
supports of $b$ and $h$ are contained in $X\times \rr d$.
\end{enumerate}

\end{prop}

\par

\begin{prop}\label{echarequiv}
Let $r,\rho \in [0,1]$ be such that $r >0$, $\omega _0\in
\mathscr P_{r,\rho}(\rr {2d})$, $a\in \SG ^{(\omega _0)}_{r,\rho }(\rr
{2d})$, and let $(x_0,\xi _0)\in (\rr d \setminus 0)\times \rr
d$. Then the following conditions are equivalent:

\medspace

\begin{enumerate}

\item $(x_0,\xi _0)\notin \Char _{(\omega _0)}^{e} (a)$;

\vrum

\item there is an element $c\in \SG^{0,0}_{r,\rho}$ which is
$e$-invertible at $(x_0,\xi _0)$, and an element $b\in \SG^{(1/\omega
_0)}_{r,\rho}$ such that $ab=c$;

\vrum

\item there is an element $c\in \SG^{0,0}_{r,\rho}$ which is
$e$-invertible at $(x_0,\xi _0)$, and elements $h\in
\SG^{-r,0}_{r,\rho}$ and $b\in \SG^{(1/\omega _0)}_{r,\rho}$ such that
\eqref{locinvop} holds;

\vrum

\item for each open conical neighbourhood $\Gamma$ of $x_0$ and open
neighbourhood $X$ of $\xi _0$, there is an element $c=\psi \otimes
\fy$ where $\psi \in \mathscr C_{x_0}^\dir(\Gamma)$ and $\fy \in
\mathscr C_{\xi _0} (X)$, and elements $h\in \mathscr S$ and $b\in
\SG^{(1/\omega _0)}_{r,\rho}$ such that \eqref{locinvop}
holds. Furthermore, the supports of $b$ and $h$ are contained in
$\Gamma \times \rr d$.

\end{enumerate}

\end{prop}

\par

\begin{prop}\label{psiecharequiv}
Let $r,\rho \in (0,1]$, $\omega _0\in \mathscr P_{r,\rho}(\rr
{2d})$, $a\in \SG ^{(\omega _0)}_{r,\rho }(\rr {2d})$, and let
$(x_0,\xi _0)\in (\rr d \setminus 0)\times (\rr d\setminus 0)$. Then
the following conditions are equivalent:

\medspace

\begin{enumerate}

\item $(x_0,\xi _0)\notin \Char _{(\omega _0)}^{\psi e} (a)$;

\vrum

\item there is an element $c\in \SG^{0,0}_{r,\rho}$ which is $\psi
e$-invertible at $(x_0,\xi _0)$, and an element $b\in \SG^{(1/\omega
_0)}_{r,\rho}$ such that $ab=c$;

\vrum

\item there is an element $c\in \SG^{0,0}_{r,\rho}$ which is $\psi
e$-invertible at $(x_0,\xi _0)$, and elements $h\in
\SG^{-r,0}_{r,\rho} + \SG^{0,-\rho}_{r,\rho}$ and $b\in \SG^{(1/\omega
_0)}_{r,\rho}$ such that \eqref{locinvop} holds;

\vrum

\item for each open conical neighbourhoods $\Gamma _1$ of $x_0$ and
$\Gamma _2$ of $\xi _0$, there is an element $c=\psi _1 \otimes \psi
_2$ where $\psi _1 \in \mathscr C_{x_0}^\dir(\Gamma _1)$ and $\psi _2
\in \mathscr C_{\xi _0}^\dir (\Gamma _2)$, and elements
$h\in \mathscr S$ and $b\in \SG^{(1/\omega _0)}_{r,\rho}$ such that
\eqref{locinvop} holds. Furthermore, the supports of $b$ and $h$ are
contained in $\Gamma _1\times \rr d$.
\end{enumerate}

\end{prop}

\par

Propositions \ref{psicharequiv} and \ref{echarequiv} follow by
the same arguments
as in the proof of Proposition 2.3 in \cite{CJT1}. Proposition
\ref{psiecharequiv} follows by similar arguments. For completeness
we give a proof of Proposition \ref{psiecharequiv} in Appendix A.

\par

As a consequence of Propositions
\ref{psicharequiv}--\ref{psiecharequiv}, we can show that the sets of
characteristic points are invariant under the choice of
pseudo-differential calculus.

\par

\begin{prop}\label{invchar}
Let $r,\rho \in [0,1]$, $\omega _0\in \mathscr P_{r,\rho}(\rr
{2d})$, $s,t\in \mathbf R$ and that $a,b\in \SG ^{(\omega
_0)}_{r,\rho}(\rr {2d})$ satisfy
$$
\op _s(a)=\op _t(b).
$$
Then the following is true:

\par

\begin{enumerate}
\item if in addition $\rho >0$, then $\Char ^\psi _{(\omega _0)}(a)
=\Char ^\psi _{(\omega _0)}(b)$;

\vrum

\item if in addition $r>0$, then $\Char ^e _{(\omega _0)}(a)
=\Char ^e _{(\omega _0)}(b)$;

\vrum

\item if in addition $r,\rho >0$, then $\Char ^{\psi e} _{(\omega
_0)}(a) =\Char ^{\psi e} _{(\omega _0)}(b)$.
\end{enumerate}

\par

\end{prop}

\par

\begin{proof}
We may assume that $s=0$, and prove only (3). The other assertions
follow by similar arguments and are left for the reader. By
\cite[Proposition 18.5.10]{Ho1}, it
follows that $b=a+h$, where $h\in \SG ^{(\omega _0/\sigma
_{r,\rho})}_{r,\rho}$ and $\sigma
_{r,\rho}$ is defined by \eqref{sigmammy}. Then for each $\ep >0$ there is a constant $R>0$ such
that $|h(x,\xi )|\le \ep \omega _0(x,\xi )$ when $|x|\ge R$ or $|\xi
|\ge R$. This implies that (3) in
Definition \ref{defchar} is fulfilled for $a$, if and only if it is
fulfilled for $b$. This gives the result.
\end{proof}

\par

\begin{rem}\label{compchar}
Let $\omega _0(x,\xi )=\eabs \xi ^r$, $r\in \mathbf R$, and assume that
$a\in \SG ^{r,0}_{1,0}(\rr {2d})$ $=\SG _{1,0}^{(\omega _0)}(\rr {2d})$
is polyhomogeneous with principal symbol $a_r\in \SG ^{r,0}_{1,0}(\rr
{2d})$ (cf. Definition 18.1.5 in \cite{Ho1}). Also let $\Char '(a)$
be the set of characteristic points of $\op (a)$ in the classical
sense (i.e., in the sense of Definition 18.1.25 in
\cite{Ho1}). Then
\begin{equation}\label{charcharrel}
\Char _{(\omega _0)}^\psi (a)\subseteq \Char '(a),
\end{equation}
where strict inclusion might appear in view of Remark 1.4 and Example
3.9 in \cite{PTT1}.

\par

By similar arguments it follows that the sets $\Char _{(\omega _0)}^e (a)$ and
$\Char _{(\omega _0)}^{\psi e} (a)$ are contained in corresponding sets of
characteristic points in \cite{CoMa}.
\end{rem}

\par

\section{Global wave-front sets}\label{sec2}

\par

In this section we define global wave-front sets for temperate distributions
with respect to Banach or Fr\'echet spaces $\cB$ and show some of their properties.
The basic ideas behind these definitions can be found in \cite{CoMa}. 

\par

We start with the following definition.

\par

\begin{defn}\label{def:wfsMB}
Let $\cB$ be a Banach or Fr\'echet space such that
$\mathscr S(\rr d)\subseteq \cB \subset \mathscr S^\prime(\rr d)$
and let $f\in \mathscr S'(\rr d)$.
\begin{enumerate}
\item The \emph{$\psi$-type} wave-front set $\WFgpB (f)$ of $f$
with respect to $\cB$ consists of all
$(x_0,\xi_0)\in \rr d
\times (\rr d\setminus 0)$ such that for every $\fy \in \mathscr
C_{x_0}(\rr d)$ and every $\psi \in \mathscr
C^\dir _{\xi _0}(\rr d\setminus 0)$ it holds
\begin{equation}\label{psiwfcondB}
\fy \cdot \psi (D)f\notin \cB\text ;
\end{equation}

\vrum

\item The \emph{$e$-type} wave-front set $\WFgeB (f)$ of $f$ with
respect to $\cB$ consists of all $(x_0,\xi_0)\in
(\rr d\setminus 0)
\times \rr d$ such that for every $\psi \in \mathscr
C^\dir _{x_0}(\rr d\setminus 0)$ and every $\fy \in \mathscr
C _{\xi _0}(\rr d)$ it holds
\begin{equation}\label{ewfcondB}
\psi \cdot \fy (D)f\notin \cB\text ;
\end{equation}

\vrum

\item The \emph{$\psi e$-type} wave-front set $\WFgpeB (f)$ of $f$
with respect to $\cB$ consists of all
$(x_0,\xi_0)\in (\rr d\setminus 0)
\times (\rr d\setminus 0)$ such that for every $\psi _1\in \mathscr
C^\dir _{x_0}(\rr d\setminus 0)$ and every $\psi _2 \in \mathscr
C ^\dir _{\xi _0}(\rr d\setminus 0)$ it holds
\begin{equation}\label{psiewfcondB}
\psi _1\cdot \psi _2 (D)f\notin \cB.
\end{equation}
\end{enumerate}

\par

\noindent
Finally, the \emph{global wave-front set} $\WFgB(f)\subseteq (\rr
d\times \rr d)\back 0$ is the set
\begin{equation*}
\WFgB(f) \equiv \WFgpB (f) \bigcup \WFgeB (f) \bigcup \WFgpeB (f).
\end{equation*}
\end{defn}

\par

From now on we assume that $\cB$ in Definition \ref{def:wfsMB} is $\SG$-admissible
(see Definition \ref{admspacesdef}). We recall that Sobolev spaces of Hilbert types or,
more generally, modulation spaces, and $\mathscr S(\rr d)$ are $\SG$-admissible.
Here we recall that the modulation spaces also contain spaces of the form
$H^2_{s_1,s_2}(\rr d)$ in \cite{CoMa,Me:94}.

\par

\begin{rem}\label{rem:order}
In a similar way to Remark 2.3 in \cite{CoMa}, we note that
Definition \ref{def:wfsMB} does not change if the conditions
\begin{align*}
\fy \cdot \psi (D)f\notin \cB,\quad \psi \cdot
\fy (D)f\notin \cB,\quad \text{and}
\quad\psi _1\cdot \psi _2 (D)f\notin \cB
\end{align*}
 in \eqref{psiwfcondB}--\eqref{psiewfcondB} are replaced by
 \begin{align*}
\psi (D)(\fy \cdot f)\notin \cB,\quad \fy
(D)(\psi \cdot f)\notin \cB,\quad \text{and}
\quad\psi_2 (D)(\psi _1\cdot f)\notin \cB,
\end{align*}
respectively (when $\cB$ is $\SG$-admissible).
In fact, let $c(x,\xi )=\fy (x)\psi (\xi )$ where
$\fy\in \mathscr C_{x_0}(\rr d)$ and $\psi \in
\mathscr C^\dir _{\xi _0}(\rr d\setminus 0)$, and let $c _1 \in \SG
^{0,0}_{r,\rho}$ be equal to $1$ on $\supp c$. Then it follows
from the symbolic calculus  that 
\begin{equation}
\label{eq:order}
\hspace*{-8pt}
\op (c_1)\op(c)=\op (c)\op (c _1)\hspace{-3mm}\mod\hspace{-1mm} \op
(\mathscr S)=\op (c)\hspace{-3mm}\mod\hspace{-1mm} \op (\mathscr S).
\end{equation}
A combination of
\eqref{eq:order} and the facts that each pseudo-differential operator
with symbol in $\SG^{0,0}_{r,\rho}$ is continuous on $\cB$ now shows that (1) in Definition
\ref{def:wfsMB} does not depend on the order we apply the operators. Here we have also
used the fact that elements in $\op (\mathscr S(\rr {2d}))$ maps $\mathscr S'(\rr d)$
continuously into $\mathscr S(\rr d) \subseteq \cB$.
\end{rem}

\par

\begin{rem}\label{remothercuttoff}
Let $f\in \mathscr S'(\rr d)$, $x_0,\xi _0\in \rr d\back 0$, $\psi
_{j,1}\in \mathscr C_{x_0}^{\dir}(\rr d\back 0)$, $\psi _{j,2}\in
\mathscr C_{\xi _0}^{\dir}(\rr d\back 0)$ for $j=1,2$ be such that
$\psi _{1,k}=1$ on $\supp \psi _{2,k}$ for $k=1,2$. Also let $\cB$ be $\SG$-admissible
with respect to $r,\rho \in [0,1]$ and $d$. If $\psi _{1,1}\cdot \psi _{1,2}(D)f\in
\cB$, then $\psi _{2,1}\cdot \psi _{2,2}(D)f\in \cB$.

\par

In fact, if $c_j=\psi _{j,1}\otimes \psi _{j,2}$, then $c_1=1$ on
$\supp c_2$, and it follows from the symbolic calculus that for some
$h\in \mathscr S$ we have
$$
\psi _{2,1}\cdot \psi _{2,2}(D)f =\op (c_2)f =\op (c_2)\op (c_1)f+\op
(h)f.
$$
The assertion now follows from the fact that $\op (c_2)$ and $\op (h)$
are continuous on $\cB$.
\end{rem}

\par

The next proposition gives an alternative
definition of the global wave-front set
in terms of intersection of sets of characteristic points described in
Section \ref{sec1}. Here we recall Definition \ref{admspacesdef} for definition
of $\SG$-ordered pairs.

\par

\begin{prop}\label{prop:charset}
Let $r,\rho \in [0,1]$, $\omega
_0\in  \mathscr P_{r,\rho}(\rr {2d})$, $(\cB ,\cC )$ be a $\SG$-ordered pair
with respect to $\omega _0$, $f\in \mathscr S'(\rr d)$, and
let
$$
\Omega = \sets {a\in \SG ^{(\omega _0)}_{r,\rho }(\rr {2d})}{\op(a)f\in \cC } .
$$
Then
\begin{alignat}{3}
\WFgpB (f) &= \bigcap _{a\in \Omega} \Char ^\psi _{(\omega _0)}(a),&
\quad &\text{when} & \quad \rho &>0\label{WFpident1}
\\[1ex]
\WFgeB (f) &= \bigcap _{a\in \Omega} \Char ^e _{(\omega _0)}(a),& \quad
&\text{when} & \quad r &>0\label{WFeident1}
\\[1ex]
\WFgpeB (f) &= \bigcap _{a\in \Omega} \Char ^{\psi e} _{(\omega
_0)}(a)& \quad &\text{when} & \quad r,\rho &>0.\label{WFpeident1}
\end{alignat}
\end{prop}

\par

\begin{proof}
Let $b\in \SG ^{(1/\omega _0)}_{r,\rho }$ and $b_1\in \SG ^{(\omega
_0)}_{r,\rho }$ be chosen such that Proposition \ref{psicontprop} (4) is fulfilled
after $a$ is replaced by $b_1$. By Remark \ref{remSGadm} it follows that
$\op _t(b_1)$ is continuous and bijective from $\cB$ onto $\cC$, with
inverse $\op _t(b)$. Since
$$
\op (b)\op (a)\in \op (\SG
^{(1)}_{r,\rho})=\op (\SG ^{0,0}_{r,\rho})
$$
when $a\in \SG^{(\omega
_0)}_{r,\rho}$, we may assume that $\omega _0=1$ and $\cC =\cB$.

\par

In order to prove \eqref{WFpident1}, we first assume that
$(x_0,\xi_0)\notin \WFgpB (f)$.  By Definition \ref{def:wfsMB}
there exist $\fy \in \mathscr C_{x_0}$ and $\psi \in \mathscr
C_{\xi _0}^\dir$ such that $a(x,\xi) \equiv \varphi (x)\psi
(\xi)\in\SG_{1,1}^{0,0}$ and $\op(a)f\in \cB$. Since
\eqref{invcond} is fulfilled with $\omega _0=1$, for some constant
$c>0$, it follows that $a$ is $\psi$-invertible at $(x_0,\xi _0)$.
Hence $(x_0,\xi _0)\notin \Char ^\psi _{(\omega _0)}(a)$, and we
have proved that $\bigcap \Char ^\psi _{(\omega _0)}(a)\subseteq
\WFgpB (f)$.

\par

It remains to prove the opposite inclusion. Let $a\in\SG^{0,0}_{r,\rho}$
be such that $(x_0,\xi_0)\notin \Char ^\psi (a)$ and
$\op(a)f\in \cB$. By Proposition \ref{psicharequiv}, there are
$\fy \in \mathscr C_{x_0}$, $\psi \in \mathscr C_{\xi _0}^\dir$,
$b\in \SG_{r,\rho}^{0,0}$ and $h\in \mathscr S$ such that
$$
\op (\fy \otimes \psi ) = \op (b)\op (a)+\op (h).
$$
Since $\op
(a)f\in \cB$, $\op (b)$ is continuous on $\cB$, and $\op (h)$
maps $\mathscr S'$ into $\mathscr S$, it follows that $\fy \cdot
(\psi (D)f)=\op (\fy \otimes \psi )f\in \cB$. Hence $(x_0,\xi
_0)\notin \WFgpB (f)$. This proves \eqref{WFpident1}. By similar
arguments we also get \eqref{WFeident1} and \eqref{WFpeident1}.
The details are left for the reader, and the proof is complete.
\end{proof}

\par

The next result describes the relation between ``regularity in $\cB$\,''
of temperate distributions and global wave-front sets:

\par

\begin{thm}\label{mainthm4}
Let $\cB$ be $\SG$-admissible, and let $f\in\cS^\prime (\rr d)$. Then
\begin{equation*}
f\in\cB \quad \Longleftrightarrow \quad
\WFgB(f)=\emptyset.
\end{equation*}
\end{thm}

\par

For the proof we need the following lemma:

\par

\begin{lemma}\label{mainthm4lemma}
Let $\cB$ be $\SG$-admissible. Then the following is true:
\begin{enumerate}
\item if $\WFgpB (f)=\emptyset$, then for each bounded open set
$X\subseteq \rr d$, there exists a non-negative
$a\in \SG ^{0,0}_{1,1}$ such that $a\ge 1$ on $X\times \rr
d$ and $\op (a)f\in \cB$;

\vrum

\item if $\WFgeB (f)=\emptyset$, then for each bounded open set
$X\subseteq \rr d$, there exists a non-negative
$a\in \SG ^{0,0}_{1,1}$ such that $a\ge 1$ on $\rr d
\times X$ and $\op (a)f\in \cB$;

\vrum

\item if $\WFgpeB (f)=\emptyset$, then for some bounded open sets
$X_1,X_2\subseteq \rr d$ such that $0\in X_1$ and $0\in X_2$, there
exists a non-negative $a\in \SG ^{0,0}_{1,1}$ such that $a\ge 1$ on
$(\rr d\setminus X_1)\times (\rr d \setminus X_2)$ and $\op (a)f\in\cB$.
\end{enumerate}
\end{lemma}

\par

\begin{proof}
We only prove (1). The other assertions follow by similar arguments
and are left for the reader.

\par

The condition $\WFgpB (f)=\emptyset$ implies that for each $(x,\xi
)\in \rr d\times (\rr d\setminus 0)$, there are functions $\fy
_{x,1}\in \mathscr C_{x}$ and $\psi _{\xi ,1}\in \mathscr C_{\xi
}^\dir$ such that $\fy _{x,1}\cdot \psi _{\xi ,1}(D)f\in\cB$.
Now let $x_0\in \rr d$ be fixed, and recall that each closed cone in
$\rr d\setminus 0$ corresponds to a compact set on the unit sphere. Hence,
by compactness, it follows that for some $\fy _{x_0,2}\in \mathscr
C_{x_0}$, $\xi _1,\dots ,\xi _N\in \rr d\back 0$ and some constant $R>0$, we have
$$
\psi _1(\xi ) = \sum _{j=1}^N \psi _{\xi _j,1}(\xi ) \ge 1,\quad
\text{when}\quad |\xi |\ge R,
$$
$\fy _{x_0,2}\otimes \psi _1\in \SG ^{0,0}_{1,1}$ and $\fy
_{x_0,2}\cdot \psi _1(D)f\in \cB$.

\par

Now choose non-negative $\fy _3\in C _0^\infty(\rr d)$ such that $\fy
_3(\xi )=1$ when $|\xi |\le R$. Then $\fy _{x_0,2}\cdot \fy
_3(D)f\in C_0^\infty(\rr d) \subseteq \cB $, since
$\fy _{x_0,2}\otimes \fy _3\in C_0^\infty(\rr {2d})$. Hence, for some
$\fy _{x_0}\in \mathscr C_{x_0}$, open
neighbourhood $U=U_{x_0}$ of $x_0$ and some constant $C>0$, the
element $a_{x_0} = C\fy _{x _0} \otimes (\psi _1+\fy _3)$ belongs
to $\SG ^{0,0}_{1,1}$ and is larger than $1$ on $U\times \rr d$.
Furthermore, $\op (a_{x_0})f\in \cB$. Summing up we have proved that
for each $x\in \rr d$, there is an open neighbourhood $U_x$ of $x$ and
an element $a_x\in \SG ^{0,0}_{1,1}$ such that $a_x\ge 1$ on $U_x$
and $\op (a_x)f\in \cB$.

\par

For each compact set $K$ we may find finite numbers of $U_{x_1},\dots
U_{x_N}$ which cover $K$. The result now follows if we choose
\[
a=a_{x_1}+\cdots +a_{x_N}.
\qedhere
\]
\end{proof}

\par

\begin{proof}[Proof of Theorem \ref{mainthm4}]
 The right implication is obvious by Definition
\ref{def:wfsMB}, since operators in $\op(\SG^{0,0}_{r,\rho})$ are
continuous from $\cB$ to itself.

\par

Assume that $\WFgB(f)=\emptyset$. Then
$\WFgpB (f)=\WFgeB (f)=\WFgpeB (f)=\emptyset$. By Lemma
\ref{mainthm4lemma} (3), there is a non-negative element $a_3\in
\SG ^{0,0}_{1,1}$, bounded open sets $X_1,X_2$ such that $0\in X_1$, $0\in
X_2$, $a_3\ge 1$ in $(\rr d\setminus X_1)\times (\rr d\setminus
X_2)$ and $\op (a_3)f\in \cB$. Furthermore, by
Lemma \ref{mainthm4lemma} (1) and (2), there are non-negative
elements $a_1,a_2\in \SG ^{0,0}_{1,1}$ such that $a_1\ge 1$ in
$X_1\times \rr d$, $a_2\ge 1$ in $\rr d\times X_2$, $\op (a_1)f\in
\cB$ and $\op (a_2)f\in \cB$. Hence, if $a=a_1+a_2+a_3$, it follows that
\begin{equation}\label{aopa}
a\in \SG ^{0,0}_{1,1},\quad
\op (a)f\in \cB\quad \text{and}\quad  a\ge
1.
\end{equation}
In particular, $a$ is elliptic in $\SG^{0,0}_{1,1}$, which implies that for some $b\in \SG
^{0,0}_{1,1}$ and $h\in \SG^{-\infty,-\infty}_{1,1}=\mathscr S$ we have
$$
\op (b)\op (a) = \operatorname{Id}+ \op (h)
$$
(cf. the proof of Proposition \ref{psiecharequiv} in Appendix A). Since $\op (b)$
and $\op (h)$ are continuous on $\cB$ and $\mathscr S \subseteq \cB$, \eqref{aopa} gives
$$
f=\op (b)\op (a)f-\op (h)f\in \cB,
$$
and the assertion follows. The proof is
complete.
\end{proof}

\par

We conclude the section by giving some remarks on
wave-front sets of modulation space type. We start to consider mapping properties
under Fourier transformation. Here it is convenient to let
$\omega _T$ be the composition of the weight $\omega \in \mathscr
P(\rr {2d})$ with the torsion $T(x,\xi )=(-\xi ,x)$, and $\mathscr
B_T$ denote the space of the pull-backs of the elements of the
translation invariant BF-space $\mathscr B$ on $\rr {2d}$ with respect to
$T$. That is,
\begin{equation}\label{torsion}
\begin{alignedat}{2}
\mathscr B_T=\sets {F\circ T}{F\in \mathscr B},\quad &\text{and}
&\quad \omega _T &=\omega \circ T ,
\\[1ex]
&\text{where} & \quad T(x,\xi ) &= (-\xi ,x).
\end{alignedat}
\end{equation}

\par
The first two equalities in the following proposition are related to
Lemma 2.4 in \cite{CoMa}.

\par

\begin{prop}
Let $\mathscr B$ be a translation invariant BF-space on $\rr {2d}$,
$\omega \in \mathscr P(\rr {2d})$, and let $T$, $\mathscr B_T$ and
$\omega _T$ be as in \eqref{torsion}. If $f\in \mathscr S'(\rr {2d})$,
then
\begin{align*}
T \big (\WF{\psi}_{M(\omega ,\mathscr B)}(f) \big ) = \WF{e}_{M(\omega _T,\mathscr
B_T)}(\widehat f),
\\[1ex]
T \big (\WF{e}_{M(\omega ,\mathscr B)}(f) \big ) = \WF{\psi}_{M(\omega _T,\mathscr
B_T)}(\widehat f),
\\[1ex]
T \big (\WF{\psi e}_{M(\omega ,\mathscr B)}(f) \big ) = \WF{\psi e}_{M(\omega
_T,\mathscr B_T)}(\widehat f).
\end{align*}
\end{prop}

\par

\begin{proof}
By Fourier's inversion formula we have
$$
|(V_\phi f)\circ T | = |V_{\widehat \phi}\widehat f|,\quad \mathscr
F(a \cdot (b(D)f)) =\check a(D)(b\cdot \widehat f).
$$
The result is now a straightforward consequence of these identities,
Remark \ref{rem:order} and the definitions. The details are left for
the reader.
\end{proof}

\medspace

Next we consider wave-front sets with respect to Fourier BF-spaces, and make
comparisons with wave-front sets of modulation space types.  In fact, in Definition
\ref{def:wfsMB} we may choose $\cB$ as the Fourier BF-space $\FB _0
{(\omega )}$, where $\mathscr B$ is a translation invariant
BF-space on $\rr d$ and $\omega \in \mathscr P(\rr{2d})$. We remark
that if $\mathscr B$ is a translation invariant BF-space on $\rr {2d}$,
then for the wave-front sets of $\psi$-type we have
\begin{equation}\label{WFeq1}
\WFpB (f)=\WF{\psi}_{\FB_0(\omega)} (f),
\end{equation}
when $\mathscr B_0$ is defined by \eqref{B0def}. This is an immediate
consequence of \eqref{compsupps}.

\par

The first type of wave-front sets with respect to general modulation
space and Fourier BF-spaces
were introduced in \cite {CJT1}. Here we recall these definitions and
show that they agree with corresponding wave-front sets of
$\psi$-type. Let $f\in \mathscr S'(\rr d)$, $\phi \in C_0^\infty (\rr
d)$ and
$\omega \in \mathscr P(\rr {2d})$. Also let $\chi_{\Gamma}$ be the
characteristic function of $\Gamma$. Then $\WFM(f)$ (denoted by
$\WFMB(f)$ in \cite{CJT1}) consists of all
pairs  $(x_0,\xi _0)\in \rr d\times (\rr d\setminus 0)$ such that
$$
\nm{(V_\phi (\fy f)) \cdot (1\otimes \chi_\Gamma )\cdot
\omega}{\mathscr B}=+\infty
$$
for every choice of open conical neighbourhood $\Gamma $
of $\xi _0$ and $\fy \in \mathscr C_{x_0}$. The wave-front set
$\WFB(f)$ (denoted by
$\WF{}_{\FB (\omega )}(f)$ in \cite{CJT1}) consists of all pairs
$(x_0,\xi _0)\in \rr d\times (\rr d\setminus 0)$ such that $|\fy
f| _{\FB(\omega,\Gamma )}=+\infty$ for every choice of open
conical neighbourhood $\Gamma $ of $\xi _0$ and $\fy \in \mathscr
C_{x_0}$. Here 
$$ 
|f|_{\FB(\omega ,\Gamma )}\equiv \nm{ \widehat
f\omega\chi _\Gamma}{\mathscr B}. 
$$

\par

\begin{prop}\label{WFidentities}
Let $f\in \mathscr S'(\rr d)$, $\mathscr B$ be a translation
invariant BF-space, $\mathscr B_0$ be defined by \eqref{B0def} and
let $\omega \in \mathscr P(\rr {2d})$. Then
$$
\WFpB (f) =
\WFM(f) = \WF{\psi}_{\FB
_0(\omega)} (f) = \operatorname{WF}'_{\FB_0(\omega )}(f).
$$
\end{prop}

\par

\begin{proof}
By Theorem 6.9 in \cite{CJT1} we have
$\operatorname{WF}'_{M(\omega,\mathscr
B)}(f)=\operatorname{WF}'_{\FB_0(\omega )}(f)$. Hence, in view of
\eqref{WFeq1}, it suffices to prove $\WF {\psi} _{\FB_0(\omega)} (f)
\!=\! \operatorname{WF}'_{\FB_0(\omega )}(f)$.
By Remark \ref{rem:order} we have
\begin{alignat*}{3}
& & (x_0,\xi _0) \notin \operatorname{WF}' _{\FB_0(\omega)} (f) &
& 
\\[1ex]
&\Longleftrightarrow &\quad
|\fy_{x_0}f|_{\FB_0(\omega,\Gamma)} &<\infty &\quad &\text{for
some}\ \fy_{x_0} \in \mathscr C_{x_0}\ \text{and} \ \Gamma =
\Gamma _{\xi _0}
\\[1ex]
&\Longleftrightarrow &\quad \nm{\psi _{\xi _0}(D)(\fy
_{x_0}f)}{\FB_0(\omega)} &<\infty & \quad &\text{for some}\ \fy
_{x_0} \in \mathscr C_{x_0}\ \text{and} \  \psi _{\xi _0} \in
\mathscr C_{\xi _0}^\dir
\\[1ex]
&\Longleftrightarrow &\quad (x_0,\xi _0) \notin \WF{\psi}_{\FB
_0(\omega)} (f). & &
\end{alignat*}
This proves the result.
\end{proof}

\par

\section{Wave-front sets for pseudo-differential
operators with smooth symbols}\label{sec3}

\par

In this section we consider mapping properties for
pseudo-differential operators with respect to global wave-front
sets. More precisely, we prove that microlocality and microellipticity
hold for pseudo-differential operators in
$\op(\SG^{(\omega_0)}_{r,\rho})$.
We start with the following result:

\par

\begin{thm}\label{mainthm2}
Let $r,\rho \in [0,1]$, $t\in \mathbf R$, $\omega _0\in \mathscr P_{r,\rho } (\rr {2d})$,
$a\in \SG^{(\omega _0)}_{r,\rho} (\rr {2d})$ and let $f\in \cS '(\rr
d)$. Moreover, let $(\cB , \cC )$ be a $\SG$-ordered pair with respect to $\omega _0$.
Then the following is true:
\begin{enumerate}
\item  if in addition $\rho >0$, then
\begin{multline}\label{pwavefrontemb1}
    \WF{\psi }_{\mathcal C} ( \op _t(a)f) \subseteq
    \WF{\psi }_{\mathcal B} (f)
          \subseteq
          \WF{\psi }_{\mathcal C} ( \op _t(a)f) \bigcup
\Char^{\psi }_{(\omega _0)}(a);
\end{multline}

\vrum

\item  if in addition $r >0$, then
\begin{multline}\label{ewavefrontemb1}
    \WF{ e}_{\mathcal C} ( \op _t(a)f) \subseteq
    \WF{ e}_{\mathcal B} (f)
          \subseteq
          \WF{e }_{\mathcal C} ( \op _t(a)f)\bigcup
\Char ^{e}_{(\omega _0)}(a);
\end{multline}

\vrum

\item  if in addition $r,\rho >0$, then
\begin{multline}\label{pewavefrontemb1}
    \WF{\psi e}_{\mathcal C} ( \op _t(a)f) \subseteq
    \WF{\psi e}_{\mathcal B} (f)
          \subseteq
          \WF{\psi e}_{\mathcal C} ( \op _t(a)f) \bigcup
\Char ^{\psi e}_{(\omega _0)}(a).
\end{multline}
\end{enumerate}
\end{thm}

\par

\begin{proof}
We only prove (3). The assertions (1) and (2) follows by similar arguments and are left for the reader.

\par

Assume that $(x_0,\xi _0)\notin \WFgpeB(f)$. We shall prove that $(x_0,\xi _0)\notin \WF{\psi
e}_{\cC} (\op (a)f)$.
For some $\psi _{1,1}\in \mathscr C_{x_0}^\dir$ and $\psi
_{1,2} \in \mathscr C_{\xi _0}^\dir$ we have
\begin{equation}\label{psiecond2}
\psi _{1,1} \cdot \psi _{1,2}(D)f\in \cB,
\end{equation}
in view of (3) in Definition \ref{def:wfsMB}. Let $\psi
_{2,1}\in \mathscr C_{x_0}^\dir$ and $\psi _{2,2} \in \mathscr C_{\xi
_0}^\dir$ be such that $\psi _{1,j}=1$ on $\supp \psi _{2,j}$, and set
$$
c_1(x,\xi ) = \psi _{1,1}(x)\psi _{1,2}(\xi )\quad \text{and}\quad
c_2(x,\xi ) = \psi _{2,1}(x)\psi _{2,2}(\xi ).
$$
Then $c_1,c_2\in \SG ^{0,0}_{r,\rho}$, and since $c_1=1$ on $\supp
c_2$, and $\SG ^{-\infty ,-\infty }_{r,\rho }=\mathscr S$, it
follows from the symbolic  calculus that
\begin{equation}\label{c2ac1}
\op (c_2)\op (a) = \op (c_2)\op (a)\op (c_1)\mod \op (\mathscr S).
\end{equation}

\par

Now we recall that the mappings
\begin{equation}\label{contmaps3}
\op (a)\, : \, \cB \to \cC,
\qquad
\op (c_2)\, : \, \cC \to \cC
\end{equation}
are continuous (cf. Proposition \ref{psicontprop}). A combination of
\eqref{psiecond2}--\eqref{contmaps3}
with the facts that $\op (c_1) =  \psi _{1,1}\cdot\psi _{1,2}(D)$
and that $\op (h)$ maps $\mathscr S'$ into $\mathscr S$ give
\begin{align*}
\psi _{2,1}\cdot \psi _{2,2}(D)(\op (a) f ) = \op (c_2)\op (a) f&
\\[1ex]
=\op (c_2)\op (a)\op (c_1)f&\mod \mathscr S \in \cC.
\end{align*}
This proves that $(x_0,\xi _0)\notin \WF{\psi e}_{\mathcal C}  (\op
(a)f)$, and the first inclusion in \eqref{pewavefrontemb1} follows.

\par

It remains to prove the second inclusion in
\eqref{pewavefrontemb1}. Assume that
$$
(x_0,\xi _0)\notin \WF{\psi e}_{\mathcal C}  (\op
(a)f)\bigcup \Char ^{\psi e}_{(\omega _0)}(a).
$$
By Remark \ref{remothercuttoff}, there exist $\psi _{1,1}\in \mathscr
C_{x_0}^\dir$ and $\psi _{1,2}\in \mathscr C_{\xi _0}^\dir$, $b\in \SG
^{(1/\omega _0)}_{r,\rho}$ and $h\in \mathscr S$ such that
$$
\psi _{1,1}\cdot \psi _{1,2}(D) (\op (a)f)\in \cC
$$
and \eqref{locinvop} holds for
$c=c_1\equiv \psi _{1,1}\otimes \psi _{1,2}$. We claim that
\begin{equation}\label{expansion1}
\op (c_2) = \op (c_2)\op (b)\op (c_1)\op (a)+ \op (h),
\end{equation}
for some $h\in \mathscr S$, where $c_2=\psi _{2,1}\otimes \psi
_{2,2}$, and $\psi _{2,1}\in \mathscr C_{x_0}^\dir$ and $\psi
_{2,2}\in \mathscr C_{\xi _0}^\dir$ are such that $\psi _{1,j}=1$ on
$\supp \psi _{2,j}$ for $j=1,2$.

\par

In fact, by combining \eqref{locinvop} with the
fact that $\op (c_2)\op (c_1)=\op (c_2)\mod \op (\mathscr S)$, we get
\begin{align*}
\op (c_2) = \op (c_2)&\op (b)\op (a)\mod \op (\mathscr S)
\\
&= \op (c_2)\op (c_1)\op (b)\op (a)\mod \op (\mathscr S)
\\
&= \op (c_2)\op (b)\op (c_1)\op (a)\mod \op (\mathscr S),
\end{align*}
and \eqref{expansion1} follows. Here the last equality follows from
the fact that 
$$
\op (c_2)[\op (b),\op (c_1)]\in \op (\mathscr S) \quad \text{and} \quad \op (\mathscr S) \op (a) \subseteq \op (\mathscr S) ,
$$
when $c_1=1$ on $\supp c_2$, where $[\cdo ,\cdo ]$ denotes the
commutator. A combination of Proposition \ref{psicontprop},
\eqref{expansion1} and the fact that $\op (c_1) (\op (a)f)\in\cC$
now shows that the mappings
\[
\op (b) \, : \, \cC \to \cB,
\qquad
\op (c_2)  \, : \, \cB \to \cB
\]
and
\begin{equation*}
\op (h)  \, : \, \mathscr S'\, \to \, \mathscr S
\end{equation*}
are continuous and that $\op (c_2)f\in \cB$.
Hence, we have showed that $(x_0,\xi _0)\notin \WFgpeB(f)$, and the proof is complete.
\end{proof}

\par

\begin{cor}
Let $r>0$, $f\in \cS'(\rr d)$ and $\fy \in C^\infty(\rr d)$ be
such that $\eabs x^{r|\alpha|}\partial^\alpha \fy(x)\in L^{\infty}(\rr
d)$ for every $\alpha$. Also let $\cB$ be $\SG$-admissible with respect to $r$, $1$ and $d$. Then
\begin{equation}\label{smooth}
\hspace*{-1.6mm}
\begin{aligned}
	\WF{\psi}_{\mathcal B}(\fy f)\subseteq \WF{\psi}_{\mathcal
        B}(f),&\; 
        \WF{e}_{\mathcal B}(\fy f)\subseteq
	\WF{e}_{\mathcal B}(f)
\\[1ex]
        \text{and}\quad \WF{\psi e}_{\mathcal B}(\fy& f)\subseteq
        \WF{\psi e}_{\mathcal B}(f).
\end{aligned}
\end{equation}
\end{cor}

\par

\begin{proof}
It follows from the assumptions that $a\equiv \fy \otimes 1\in
SG^{0,0}_{r,1}$. Hence Theorem \ref{mainthm2} gives
$$
\WF{\psi e}_{\mathcal B}(\fy f)=\WF{\psi e}_{\mathcal B}(\Op(a)
f)\subseteq \WF{\psi e}_{\mathcal B}(f),
$$
which is the last inclusion in \eqref{smooth}. The other inclusions
follow by similar arguments. The proof is complete.
\end{proof}

\par

Next we apply Theorem \ref{mainthm2} on
operators which are elliptic with respect to $\omega _0\in \mathscr
P_{\rho,\delta}(\rr {2d})$ when $0< r ,\rho \le 1$. We recall that $a$
and $\op (a)$ are called $\SG$-\emph{elliptic} with respect to
$\SG^{(\omega _0)}_{r,\rho}(\rr {2d})$ or $\omega _0$, if
there is a compact set $K \subset \rr{2d}$ and a positive constant $c$
such that \eqref{invcond} holds when $(x,\xi )\notin K$.
Since $|a(x,\xi )|\le C\omega _0(x,\xi )$, it follows from the
definitions that for each multi-index $\alpha$, there are constants
$C_{\alpha,\beta}$ such that
\begin{equation*}
|D ^\alpha _xD ^\beta _\xi a(x,\xi )| \le C_{\alpha
,\beta}|a(x,\xi )|\eabs x^{-r |\alpha |}\eabs \xi ^{-\rho |\beta |},
\qquad (x,\xi)\in \rr {2d}\setminus K,
\end{equation*}
when $a$ is $\SG$-elliptic (see, e.g., \cite {Ho1,BBR}).

\par

It follows from Lemma \ref{mainthm4lemma} that $\Char _{(\omega
_0)}(a)=\emptyset$ if and only if $a$ is $\SG$-elliptic with respect
to $\omega _0$. The following result is now an immediate consequence
of Theorem \ref{mainthm2}:

\par

\begin{thm}\label{hypoellthm}
Let $r,\rho \in (0,1]$, $t\in \mathbf R$, $\omega
_0\in \mathscr P_{r,\rho } (\rr {2d})$, $a\in \SG^{(\omega
_0)}_{r,\rho} (\rr d)$ be $\SG$-elliptic with respect to $\omega _0$
and let $f\in \cS '(\rr d)$. Moreover, let $(\cB ,\cC
)$ be a $\SG$-ordered pair with respect to $\omega _0$. Then
\begin{align*}
\WF{\psi} _{\mathcal C} (\op _t(a)f) &=
\WF{\psi}_{\mathcal B} (f),
\\[1ex]
\WF{e}_{\mathcal C} (\op _t(a)f) &=
\WF{e}_{\mathcal B} (f),
\\[1ex]
\WF{\psi e}_{\mathcal C} (\op _t(a)f) &= \WF{\psi
e}_{\mathcal B} (f).
\end{align*}
\end{thm}

\par

\begin{thm}\label{corhypoellthm}
Assume that the hypothesis in Theorem \ref{hypoellthm} is fulfilled,
let $g\in \cC$ and let $f\in \mathscr
S'(\rr d)$ be a solution to the equation
\begin{equation}\label{psdiffeq}
\op _t(a)f = g.
\end{equation}
Then $f\in \cB$.
\end{thm}

\par

\begin{proof}
The result follows by combining Theorems \ref{mainthm4} and
\ref{hypoellthm}.
\end{proof}

\par

\section{Examples}
\label{sec:ex}

In this section we show how the results proved above can be applied to
problems involving partial differential operators. In particular, we focus
on the sets of characteristic points and some links on possible application
in numerical analysis.

\par

In the first examples we consider wave-front properties for hypoelliptic problems,
especially for hypoelliptic partial differential operators with constant coefficients.

\par

\begin{example}\label{heatexample}
Let $a(x,\xi )$ be the symbol of a linear partial differential operator on $\rr d$ with
constant coefficients, which is hypoelliptic in the sense of \cite{Ho1}. Then
$a(x,\xi )=a_2(\xi )$ for some $a_2$. Furthermore, $a$ is $\psi$- and
$\psi e$-elliptic with respect to
\begin{equation}\label{schrweight}
\omega _0(x,\xi ) = \omega _0(\xi )=(1+|a_2(\xi )|^2)^{1/2},
\end{equation}
which belongs to $\mathscr P_{r,\rho}(\rr {2d})$, for some $r,\rho >0$. In particular we may 
apply Theorems \ref{hypoellthm} and \ref{corhypoellthm} on $\op (a)$.

\par

For the set of  characteristic point with respect to the $e$ component we have
$$
\Char^{e}_{(\omega_0)}(a)\subseteq (\rr{d}\setminus 0)\times \{(0,0)\}.
$$
Hence for every $f\in \mathscr S'(\rr d)$ and any modulation space
$M(\omega ,\mathscr B)$ we have
\begin{align*}
\WF {\psi} _{M(\omega /\omega _0,\mathscr B)}(\Op (a)f)
&=
\WF {\psi}_{M(\omega ,\mathscr B)}(f),
\\[1ex]
\WF {\psi e}_{M(\omega /\omega _0,\mathscr B)}(\Op (a)f)
&=
\WF {\psi e}_{M(\omega ,\mathscr B)}(f),
\\[1ex]
\WF {e}_{M(\omega /\omega _0,\mathscr B)}(\Op (a)f) 
&\subseteq 
\WF {e}_{M(\omega ,\mathscr B)}(f)
\\[1ex]
&\subseteq 
\WF {e}_{M(\omega /\omega _0,\mathscr B)}(\Op (a)f)\cup \big ( (\rr{d}\setminus 0)\times
\{ 0\} \big ).
\end{align*}

\par

Now we examine what kind of information for a parametrix $E$ to $\op (a)$ we may
obtain from the local wave-front property. From Proposition \ref{psicontprop} (3) it
follows that $\op (a)$ maps $M^{p,q}_{(\omega _0)}$ to $M^{p,q}$. Since
$$
M^{p,q}\cap \mathscr E'= \mathscr F\!L^q\cap\mathscr E'
$$
and $\op (a)E$ is locally in $\mathscr FL^\infty$, it follows that $E$ is locally in
$\mathscr FL^\infty _{(\omega _0)}$. This means that for each $\fy\in C_{0}^\infty$,
there is a constant $C$ such that
\begin{equation}\label{parametrixest}
|\mathscr F(\fy E)(\xi )|\le C\omega_0(\xi )^{-1}.
\end{equation}

\medspace

An interesting case appears when $x\in \rr d$ is replaced by $(x,t)\in \rr {d+1}$,
and $\op (a)$ agrees with the heat operator
$$
\partial_t-\Delta_x,\qquad (x,t)\in \rr {d+1} ,
$$
which is a classical example on a hypoelliptic operator.
The symbol is given by $a(x,t,\xi,\tau)=|\xi|^2 +i\tau$. In this case, $\omega _0$
takes the form
\begin{equation}\tag*{(\ref{schrweight})$'$}
\omega_0 =(1+|\xi|^4 +|\tau|^2)^{1/2}.
\end{equation}
Hence, \eqref{parametrixest} shows that if $E$ is a fundamental solution to
$\op (a)$, then for each  $\fy\in C_{0}^\infty$, there is a constant $C$ such that
$$
|\mathscr F(\fy E)(\xi ,\tau )|\le C(1+|\xi|^4 +|\tau|^2)^{-1/2}.
$$
This can easily be verified by numerical computations.
\end{example}

\par

Next we show how a small disturbance of any hypoelliptic operator in the
previous example makes that all the sets of characteristic points become empty.

\par

\begin{example}\label{examplheat}
Let $0<r,\rho \le 1$, $a(x,\xi )=a_1(x,\xi )+a_2(\xi )$ be such that
the following conditions are fulfilled:
\begin{enumerate}
\item $|a(x,\xi )|\ge c$ for some constant $c>0$ outside a compact set in $\rr{2d}$;

\vrum

\item $a_1\in \SG ^{0,0}_{r,\rho}(\rr {2d})$;

\vrum

\item $a_2$ is the symbol of a linear partial differential operator
with constant coefficient which is hypoelliptic in the sense of
\cite{Ho1}.
\end{enumerate}
Then $a$ is elliptic with respect to $\omega _0(x,\xi )$ in \eqref{schrweight}. Hence we may 
again apply
Theorems \ref{hypoellthm} and \ref{corhypoellthm} on $\op (a)$.

\par

An interesting case concerns the modified heat operator
$a_1(x,t)+\partial _t-\Delta _x$, where $(x,t)\in \rr {d+1}$,
$a_1(x,t)\in C^\infty (\rr {d+1})$ and $a_1$ is equal to $c>0$ outside
a compact set in $\rr {d+1}$. The symbol of the operator is $a(x,t,\xi
,\tau )=a_1(x,t) +|\xi |^2+i\tau$. In this case, $a$ is elliptic with
respect to \eqref{schrweight}$'$.
Hence, if $\omega \in \mathscr P(\rr {2d})$, $\mathscr B$ is a
translation invariant BF-space on $\rr d$, and \eqref{psdiffeq} holds
for some $f\in \mathscr S'(\rr d)$ and $g\in M(\omega /\omega
_0,\mathscr B )$, then it follows from Theorem \ref{corhypoellthm}
that $f\in M(\omega ,\mathscr B )$.
\end{example}

\par

In the next example we consider the fundamental solutions of the
Schr{\"o}dinger operator.

\par

\begin{example}\label{Schrodexample}
The Schr{\"o}dinger operator for a free particle has the form 
$$
i\partial_t-\Delta_x,\qquad (x,t)\in \rr {d+1} ,
$$
and the symbol is given by $a(x,t,\xi,\tau)=|\xi|^2 -\tau$. Let $\omega _0$ be
the same as in \eqref{schrweight}. Then, $a\in\mathscr \SG^{(\omega_0)}_{1,1}$,
and the sets of characteristic points are
\begin{align*}
&\Char^{\psi}_{(\omega_0)}(a)=\sets{(x,t,\xi,\tau)\in \rr {d+1}\times (\rr{d+1}\setminus 0)}{\xi=0, 
\tau>0},
\\[1ex]
&\Char^{e}_{(\omega_0)}(a)=\sets{(x,t,\xi,\tau)\in (\rr{d+1}\setminus 0)\times \rr {d+1}}{\tau=|\xi|^2},
\\[1ex]
&\Char^{\psi e}_{(\omega_0)}(a)=\sets{(x,t,\xi,\tau)\in (\rr{d+1}\setminus 0)\times 
(\rr{d+1}\setminus  0)}{\xi=0,\tau> 0}.
\end{align*}
This implies that numerical computations can be rather easily performed as long
as the frequency $\tau$ related to the $t$-variable is negative.

\par

Let $E$ denote the fundamental solution of $\Op (a)$. From Proposition
\ref{psicontprop} (3) it follows that $\op (a)$ maps $M^{p,q}_{(\omega _0)}$
to $M^{p,q}$. Now let $\Gamma$ be a cone which does not hit the set
$\sets{(0,\tau)}{\tau>0}$. Then the same arguments as in Example
\ref{heatexample} shows that
$$
\nm {\mathscr F(\fy E)\psi \omega_0}{L^\infty}<\infty.
$$
for every $\fy \in C_{0}^{\infty}$ and $\psi \in\mathscr C^{\dir}_{\xi_0}(\Gamma)$.
\end{example}

\par

\begin{rem}
Here we note, once again, that the sets of characteristic points $\Char'(a)$, as
defined in Section 18.1 in \cite{Ho1}, for the operators in Examples
\ref{Schrodexample} and \ref{heatexample}, are strictly larger comparing to
the corresponding local components $\Char^{\psi}_{(\omega_0)}(a)$. In fact,
if $a$ is the same as in Example \ref{Schrodexample} or \ref{heatexample}, then
$$
\Char'(a)=\sets{(x,t,\xi,\tau)\in \rr {d+1}\times (\rr{d+1}\setminus 0)}{\xi=0,\tau\neq 0},
$$
which strictly includes $\Char^{\psi}_{(\omega_0)}(a)$.
\end{rem}

\par

In the next example we consider consider an operator with polynomial coefficients, with different growth
with respect to spatial axis. We show the associated \emph{loss of decay} at infinity, described through
our wave-front sets, where the reference spaces are the weighted Sobolev spaces $H^2_{s,t}(\rr d)$.

\par

\begin{example}\label{wSobolev}
Let $a(x,\xi )$ be the symbol on $\rr 4$, given by
$$
a(x,\xi)=(1+x_1^2+x_2^4)(1+|\xi|^2),\quad x,\xi \in \rr 2.
$$
Also assume that $f,g\in \mathscr S '(\rr 2)$ are chosen such that the equation $\op (a)f=g$ is fulfilled.
By letting $\omega _0(x,\xi )$ be the weight which is equal to $a$, it follows that
$a$ is elliptic with respect to $\omega _0$, and that $a\in
\SG ^{(\omega _0)}_{1,1}(\rr 4)$. If
\begin{multline*}
H^2_{(\omega _0)}(\rr 2) \equiv \sets {f\in \mathscr S '(\rr 2)}{(1+x_1^2+x_2^4)(1-\Delta )f\in L^2(\rr 2)}
\\[1ex]
= M^2_{(\omega _0)}(\rr 2),
\end{multline*}
then it is obvious that $f\in H^2_{(\omega _0)}$, if and only if $g\in L^2$, which is also confirmed by
Theorem \ref{hypoellthm}.

\par

Furthermore, $a\in \SG ^{4,2}_{1,1}(\rr 4)$ is
$\SG$ hypoelliptic with \textit{inverse order} $(2,2)$, in the sense of \cite{Co},
since, for $|x|+|\xi|\ge R>0$ and a suitable $C>1$,
$$
C^{-1}\norm{x}^2\norm{\xi}^2\le a(x,\xi)\le C\norm{x}^4\norm{\xi}^2.
$$
We consider $\op (a)$ acting on the scale of weighted Sobolev spaces $H^2_{s,t}(\rr 2)$,
$s,t\in\mathbf R$, and the solutions of the equation $\op (a)f = g$ with
$g\in H^2_{0,0}(\rr 2)=L^2(\rr 2)$. Let 
\begin{align*}
	c_1&=\psi_1\otimes\varphi_1, \quad \text{where}\quad \psi_1 \in{\mathscr C}^\dir_{(x_1,0)},\ 
	\varphi_1\in{\mathscr C}_{\xi_0},\  x_1\not=0,\ \xi_0\in\rr 2,
\intertext{and}
	c_2&=\psi_2\otimes\varphi_2, \quad \text{where}\quad \psi_2 \in{\mathscr C}^\dir_{(0,x_2)},\ 
	\varphi_2\in{\mathscr C}_{\xi_0},\  x_2\not=0,\  \xi_0\in\rr 2.
\end{align*}
Then
there exist $\SG$ symbols $b_j\in\SG^{-2j,-2}_{1,1}(\rr 4)$ and $h_j\in\mathscr S(\rr 4)$
for $j=1,2$ such that $\op(b_j)\op(a)=\op(c_j)+\op(h_j)$. (In fact, $\op (a)$ is invertible with
inverse $\op (b)$, where $b\in \SG ^{(1/\omega _0)}_{1,1}(\rr 4)$.) This implies that for some choices
of $f$ and $g$, there are \emph{anisotropic loss of decays} for $f$, in the sense that
for any $\varepsilon >0$,
\begin{align*}
	((x_1,0),\xi)&\in\WF{e}_{H^2_{2,2+\varepsilon}}(f)\setminus \WF{e}_{H^2_{2,2}} (f),\  x_1\not=0,
	\\[1ex]
	((0,x_2),\xi)&\in\WF{e}_{H^2_{2,4+\varepsilon}}(f)\setminus \WF{e}_{H^2_{2,4}} (f),\  x_2\not=0,
\end{align*}
and the same results hold for the corresponding $\WF{\psi e}_{H^2_{s,t}}(f)$ components, by
choosing symbols
$$
c_j=\psi_{1j}\otimes\psi_{2j},\quad \text{where}\quad \psi_{1j}\in {\mathscr C}^\dir_{x},
\ \psi_{2j}\in{\mathscr C}^\dir_{\xi},\  j=1,2,\ x,\xi \not=0.
$$
In particular, the first formula holds also for any $(x_1,x_2)\in\rr 2$, $x_1\not=0$.
\end{example}

\par

In the next example, we apply the theory to the propagation of electromagnetic waves in a wave-guide. The corresponding problem, concerning electromagnetic waves in a uniform cylindrical waveguide, has been considered by Kristensson in \cite{Kr} and by Khrennikov-Nilsson in \cite {KhN}, and will be viewed here in the context of pseudo-differential operators and wave-front sets. 

\par

\begin{example}
Let $S\subset \rr{3}$ be a conducting cylindrical tube along the $x_3$ axis with $\partial S$ as its boundary surface. With $\Omega$ we
denote the cross section in the $x_1 x_2$-plane, assumed to be compact and with a smooth boundary.

\par

For the propagation of the electromagnetic wave, we then have a scalar field $\phi=\phi(x,t)=\phi(x_1,x_2,x_3,t)$,
 $(x_1,x_2,x_3)\in S$, $t\in\rr{}$, satisfying the wave equation in the waveguide, namely
$$
\Big(\frac{\partial^2 }{\partial x_1^2}+\frac{\partial^2 }{\partial x_2^2}+\frac{\partial^2 }{\partial x_3^2} -\frac{\partial^2 }{\partial t^2} \Big)\phi=0.
$$
As boundary conditions we may choose either the Neumann boundary condition, i.e. let $\partial \phi /\partial \overline n=0$ on $\partial S$, or the Dirichlet boundary condition, i.e. let $\phi=0$ on $\partial S$. To solve this type of equation one can use a splitting technique, separating the pairs of variables $(x_1,x_2)$ and $(x_3,t)$. 
The solution can then be written as 
$$
\phi(x_1,x_2,x_3,t)=\sum_{n=1}^{\infty}v_n(x_1,x_2) w_n(x_3,t),
$$
where $\{v_n(x,y)\}$, $n=1,2,\dots$, is a complete orthonormal system of eigenfunctions of 
$$
\Big(\frac{\partial^2 }{\partial x_1^2}+\frac{\partial^2 }{\partial x_2^2}+m^2_n\Big)v_n(x_1,x_2)=0,
$$
with $\partial v_n /\partial \overline n=0$ or $v_n=0$ on $\partial \Omega$,
with eigenvalues $-m_n^2$. 
We have that $w_n$ satisfies the Klein-Gordon equation
\begin{equation}\label{waveeq}
\Big(\frac{\partial^2 }{\partial x_3^2}-\frac{\partial^2 }{\partial t^2}-m^2_n\Big)w_n(x_3,t)=0.
\end{equation}

\par

In the $(x_1,x_2)$ direction the problem can be solved using discretization and Fourier series, since the cross section is finite. We therefore consider equation \eqref{waveeq}. The operator on the left-hand side of \eqref{waveeq} can be written as a pseudo-differential operator with symbol 
$a(x_3,t,\xi_3,\tau)=\tau^2-\xi_3^2-m_n^2$. Let $\omega_0(x_3,t,\xi_3,\tau)=(1+\xi_3^2+\tau^2)$. Then $a\in \mathscr \SG^{(\omega_0)}_{1,1}$ has the following sets of characteristic points 
\begin{align*}
&\Char^{\psi}_{(\omega_0)}(a)=\sets{(x_3,t,\xi_3,\tau)\in \rr 4}{\tau=\pm \xi_3\neq 0},
\\[1ex]
&\Char^{e}_{(\omega_0)}(a)=\sets{(x_3,t,\xi_3,\tau)\in (\rr{2}\setminus 0)\times \rr 2}{\tau=\pm\sqrt{\xi_3^2+m_n^2}},
\\[1ex]
&\Char^{\psi e}_{(\omega_0)}(a)=\sets{(x_3,t,\xi_3,\tau)\in (\rr{2}\setminus 0)\times \rr{2}}{\tau=\pm \xi_3\neq 0}.
\end{align*}

\par

Now let $\Gamma$ be a closed cone which does not contain any point from the set $\sets{(\xi_3,\tau)}{\tau=\pm \xi_3}.$ Furthermore let $w_n$ be the solution of equation \eqref{waveeq}. Since the right-hand side in \eqref{waveeq} is zero,
it follows that for every $x_0$ there exist $\fy \in \mathscr C_{x_0}$ and $\psi \in\mathscr C^{\dir}_{\xi_0}(\Gamma)$ such that 
$\psi\mathscr F(\fy w_n)\in \mathscr S.$ In particular, for every $x_0$ there exist $\fy\in \mathscr C_{x_0}$ and $\psi \in\mathscr C^{\dir}_{\Gamma}$ such that 
$$
\nm{\mathscr F(\fy w_n)\psi \omega_0}{L^q}<\infty.
$$
\end{example}

\par

\section{Wave-front sets with respect to sequences of spaces}\label{sec4}

\par

In this section we define wave-front sets based on sequences of
admissible spaces, and discuss basic results. In the first part we consider sequences of spaces which are parameterized with one index, and prove that the mapping properties in Section \ref{sec3} extend to wave-front sets with respect to sequences. Thereafter we discuss further extensions where we consider sequences of spaces which are parameterized with two indices. In the last part we give some examples on new types of wave-front sets that can be constructed, and show some consequences of our investigations. For example, here we introduce wave-front sets which are related to "classical wave-front sets" in the sense that they are wave-front sets with respect to classical spaces of smooth functions. In particular, we show that (a refinement of) the wave-front set of Schwartz-type
treated in \cite{CoMa} can be obtained as a wave-front set based on sequences of admissible spaces.

\par

\subsection{Wave-front sets with respect to sequences with one index parameter}

\par

We start with defining wave-front sets with respect to sequences of spaces $(\cB _j)$ depending on the parameter $j\in J$.

\par

\begin{defn}\label{defsuperposWF}
Let $J$ be an index set, $(\cB_j)\equiv(\cB_j)_{j\in J}$ be a sequence of
Banach or Fr{\'e}chet spaces such that $\mathscr S(\rr d)\subseteq \cB _j\subseteq \mathscr S'(\rr d)$ holds for every $j\in J$,  and let $f\in \cS^\prime (\rr d)$.
\begin{enumerate}
\item The \emph{$(\psi ,\cup )$-type} (\emph{$(\psi ,\cap )$-type}) wave-front set $\WF {\psi ,\, \cup}_{(\mathcal B _j)}(f)$
 ( $\WF {\psi ,\, \cap}_{(\mathcal B _j)}(f)$) 
with respect to $(\mathcal B _j)$, consists of all pairs $(x_0,\xi
_0)$ in $\rr d\times (\rr d\setminus 0)$ such that for every $\fy \in
\mathscr C_{x_0}(\rr d)$ and every $\psi \in \mathscr
C^\dir _{\xi _0}(\rr d\setminus 0)$ it holds
\begin{equation}\label{psiwfpsicond2}
\fy \cdot \psi (D)f\notin \cB_j
\end{equation}
for some $j\in J$ (for every $j\in J$);

\vrum

\item The \emph{$(e ,\cup )$-type} (\emph{$(e ,\cap )$-type}) wave-front set $\WF {e ,\, \cup}_{(\mathcal B _j)}(f)$
 ( $\WF {e ,\, \cap}_{(\mathcal B _j)}(f)$) with
respect to $(\mathcal B _j)$, consists of all pairs $(x_0,\xi _0)$ in
$(\rr d\setminus 0)\times \rr d$ such that for every $\psi \in \mathscr
C^\dir _{x_0}(\rr d\setminus 0)$ and every $\fy \in \mathscr
C_{\xi _0}(\rr d)$ it holds
\begin{equation}\label{psiwfecond2}
\psi \cdot \fy (D)f\notin \cB_j
\end{equation}
for some $j\in J$ (for every $j\in J$);

\vrum

\item The \emph{$(\psi e,\cup )$-type} (\emph{$(\psi e ,\cap )$-type}) wave-front set  $\WF {\psi e ,\, \cup}_{(\mathcal B _j)}(f)$
 ( $\WF {\psi e ,\, \cap}_{(\mathcal B _j)}(f)$)
with respect to $(\mathcal B _j)$, consists of all pairs
$(x_0,\xi _0)$ in $(\rr d\setminus 0)\times (\rr d\setminus 0)$ such
that for every $\psi_1 \in \mathscr C^\dir_{x_0}(\rr d\setminus 0)$
and every $\psi_2 \in \mathscr C^\dir _{\xi _0}(\rr d\setminus 0)$ it
holds
\begin{equation}\label{psiwfpsiecond2}
\psi _1\cdot \psi _2 (D)f\notin \cB_j
\end{equation}
for some $j\in J$ (for every $j\in J$).
\end{enumerate}

\noindent
Finally, the \emph{global wave-front sets of $\cup$- and $\cap$-types} $\WF{\cup}_{(\mathcal B
_j)}(f)\subseteq (\rr d\times \rr d)\back 0$ and $\WF{\cap}_{(\mathcal
B _j)}(f)\subseteq (\rr d\times \rr d)\back 0$ are the sets
\begin{equation*}
\WF{\cup}_{(\mathcal B _j)}(f) \equiv \WF {\psi ,\, \cup}_{(\mathcal B
_j)}(f)\bigcup \WF {e ,\, \cup}_{(\mathcal B _j)}(f) \bigcup \WF {\psi e ,\,
\cup}_{(\mathcal B _j)}(f),
\end{equation*}
and
\begin{equation*}
\WF{\cap}_{(\mathcal B _j)}(f) \equiv \WF {\psi ,\, \cap}_{(\mathcal B
_j)}(f)\bigcup \WF {e ,\, \cap}_{(\mathcal B _j)}(f) \bigcup  \WF {\psi e
,\, \cap}_{(\mathcal B _j)}(f)
\end{equation*}
respectively.
\end{defn}

\par

\begin{example}
We can consider wave-front sets with respect to sequences of the form
\begin{equation}\label{Jsetdef}
(\mathcal B _j)\equiv (\mathcal B _j) _{j \in J},\quad
\text{with}\quad \mathcal B _j = M(\omega _j,\mathscr B_j ),
\end{equation}
where $\omega _j \in \mathscr P(\rr {2d})$, $\mathscr B_j$ is a
translation invariant BF-space on $\rr d$, and $j$ belongs to some
index set $J$.
\end{example}

\par

\begin{rem}\label{remhormwfsets}
Let $p_j,q_j \in [1,\infty]$, $\mathscr B_j=L_1^{p_j,q_j}(\rr {2d})$,
$\omega _j(x,\xi ) = \eabs {x,\xi}^{-j}$ and let $\mathcal B_j$ be as
in \eqref{Jsetdef} for $j\in J=\mathbf N_0$. Then it follows that $\WF
{\psi,\, \cup}_{(\mathcal B _j)}(f)$,  $\WF {e,\, \cup}_{(\mathcal B
_j)}(f)$ and  $\WF {\psi e,\, \cup}_{(\mathcal B _j)}(f)$ in
Definition \ref{defsuperposWF} are
equal to the wave-front sets $\WF {\psi}(f)$,  $\WF {e}(f)$ and  $\WF
{\psi e}(f)$  in \cite{CoMa}, respectively. In particular, it follows
that $\WF {\cup}_{(\mathcal B _j)}(f)$ is equal to the global
wave-front set $\WF {}_{\mathscr S}(f)$, which in \cite{CoMa} is
denoted by $\WF {}_{\mathcal S}(f)$.
\end{rem}

\par

\begin{rem}\label{supinfequalwf}  
Obviously, if $\cB_j=\cB$ for
every $j\in J$, then
$$
\WF {\psi ,\, \cup}_{(\mathcal B _j)}(f)=\WF {\psi ,\,
\cap}_{(\mathcal B _j)}(f) =\WF {\psi}_{\cB}(f),
$$
and similarly for $\WF {e,\, \cup}_{(\mathcal B _j)}(f)$ and $\WF
{\psi e,\, \cup}_{(\mathcal B _j)}(f)$.
\end{rem}

\par

In the following two results we make some basic remarks, which also motivates the notations for wave-front sets of sequence types. 

\par

\begin{prop}\label{capcupprop1}
Let $\cB _j$ be the same as in Definition \ref{defsuperposWF}, and let $f\in \mathscr S'(\rr d)$. Then
\begin{equation*}
\begin{alignedat}{2}
\bigcup \WF {\psi }_{\mathcal B _j}(f) &\subseteq \WF {\psi ,\, \cup}_{(\mathcal B _j)}(f),& \qquad
\bigcup \WF {e }_{\mathcal B _j}(f) &\subseteq \WF {e ,\, \cup}_{(\mathcal B _j)}(f)
\\[1ex]
\bigcup \WF {\psi e}_{\mathcal B _j}(f) &\subseteq \WF {\psi e,\, \cup}_{(\mathcal B _j)}(f),&
\qquad
\bigcup \WF {}_{\mathcal B _j}(f) &\subseteq \WF {\cup}_{(\mathcal B _j)}(f),
\end{alignedat}
\end{equation*}
and
\begin{equation*}
\begin{alignedat}{2}
\bigcap \WF {\psi }_{\mathcal B _j}(f) &= \WF {\psi ,\, \cap}_{(\mathcal B _j)}(f),& \qquad
\bigcap \WF {e }_{\mathcal B _j}(f) &= \WF {e ,\, \cap}_{(\mathcal B _j)}(f)
\\[1ex]
\bigcap \WF {\psi e}_{\mathcal B _j}(f) &= \WF {\psi e,\, \cap}_{(\mathcal B _j)}(f),&
\qquad
\bigcap \WF {}_{\mathcal B _j}(f) &= \WF {\cap}_{(\mathcal B _j)}(f).
\end{alignedat}
\end{equation*}
\end{prop}

\par

Proposition \ref{capcupprop1} is a straightforward consequence of the definitions. The details are left for the reader. We remark that we may choose $\cB _j$ and $f\in \mathscr S'(\rr d)$ such that equality is not attained in the first inclusion in Proposition \ref{capcupprop1} (cf. \cite[Example 1.11]{PTT2}). On the other hand, the following generalization of Theorem \ref{mainthm4} shows that equality is attained when one of these sets are empty, provided $\cB _j$ in addition are $\SG$-admissible.

\par

\renewcommand{\rubrik}{Theorem \ref{mainthm4}$'$}

\begin{tom}
Let $\cB _j$ be $\SG$-admissible for every $j$, and let $f\in\cS^\prime (\rr d)$. Then
\begin{alignat*}{5}
f&\in \bigcup \cB _j &\quad &\Longleftrightarrow &\quad \bigcap \WF {}_{\cB _j}(f) &= \emptyset &\quad &\Longleftrightarrow &\quad \WF {\cap}_{(\cB _j)}(f) &= \emptyset ,
\intertext{and}
f&\in \bigcap \cB _j &\quad &\Longleftrightarrow &\quad \bigcup \WF {}_{\cB _j}(f) &= \emptyset &\quad &\Longleftrightarrow &\quad \WF {\cup}_{(\cB _j)}(f) &= \emptyset .
\end{alignat*}
\end{tom}

\par

\begin{proof}
We only prove the third and fourth equivalences. The first two equivalences follow by similar arguments and are left for the reader.

\par

The third equivalence follows immediately from the definitions and Theorem \ref{mainthm4}. In view of Proposition \ref{capcupprop1}, the result therefore follows if we prove that $f\in \bigcap \cB _j$ implies $\WF {\cup}_{(\cB _j)}(f) = \emptyset$.

\par

Therefore, assume that $f\in \bigcap \cB _j$. Since every $\cB _j$ is $\SG$-admissible, and $g_1\otimes g_2\in \SG ^{0,0}_{1,1}$ when $g_1$ and $g_2$ are any cutoff or directional cutoff, it follows that $g_1\cdot g_2(D)f\in \cB _j$ for every $j$. This implies that $\WF {\cup}_{(\cB _j)}(f) = \emptyset$, and the result follows.
\end{proof}

\medspace

\subsection{Wave-front sets with respect to sequences of spaces with two indices parameters}

\par

Next we shall consider wave-front sets with respect to sequences of spaces, parameterized with two indices. We start with defining wave-front sets with respect to sequences of spaces $(\cB _{j,k})$ depending on the parameters $j,k\in J$.

\par

\begin{defn}\label{defsuperposWF2}
Let $J$ be an index set, let $(\cB_{j,k})\equiv (\cB_{j,k})_{j,k\in J}$ be a sequence of
Banach or Fr{\'e}chet spaces such that $\mathscr S(\rr d)\subseteq \cB _{j,k}\subseteq \mathscr S'(\rr d)$ holds for every $j,k\in J$,  and let $f\in \cS^\prime (\rr d)$.
\begin{enumerate}
\item The \emph{$(\psi ,\cupap )$-type} (\emph{$(\psi ,\capup )$-type}) wave-front set $\WF {\psi ,\, \cupap}_{(\mathcal B _{j,k})}(f)$
 ( $\WF {\psi ,\, \capup}_{(\mathcal B _{j,k})}(f)$) 
with respect to $(\mathcal B _{j,k})$, consists of all pairs $(x_0,\xi
_0)$ in $\rr d\times (\rr d\setminus 0)$ such that for every $\fy \in
\mathscr C_{x_0}(\rr d)$ and  $\psi \in \mathscr
C^\dir _{\xi _0}(\rr d\setminus 0)$, there exists $j\in J$ such that (and for every $j\in J$ it holds)
\begin{equation}\label{psiwfpsicond3}
\fy \cdot \psi (D)f\notin \cB_{j,k}
\end{equation}
for every $k\in J$ (for some $k\in J$);

\vrum

\item The \emph{$(e ,\cupap )$-type} (\emph{$(e ,\capup )$-type}) wave-front set $\WF {e ,\, \cupap}_{(\mathcal B _{j,k})}(f)$ ( $\WF {e ,\, \capup}_{(\mathcal B _{j,k})}(f)$)
with respect to $(\mathcal B _{j,k})$, consists of all pairs $(x_0,\xi
_0)$ in $(\rr d\setminus 0)\times \rr d$ such that for every $\psi \in \mathscr
C^\dir _{x_0}(\rr d\setminus 0)$ and $\fy \in \mathscr
C_{\xi _0}(\rr d)$, there exists $j\in J$ such that (and for every $j\in J$ it holds)
\begin{equation}\label{psiwfecond3}
\psi \cdot \fy (D)\notin \cB_{j,k}
\end{equation}
for every $k\in J$ (for some $k\in J$);

\vrum

\item The \emph{$(\psi e ,\cupap )$-type} (\emph{$(\psi e ,\capup )$-type}) wave-front set $\WF {\psi e ,\, \cupap}_{(\mathcal B _{j,k})}(f)$
 ( $\WF {\psi e ,\, \capup}_{(\mathcal B _{j,k})}(f)$)
with respect to $(\mathcal B _{j,k})$, consists of all pairs $(x_0,\xi
_0)$ in $(\rr d \setminus 0) \times (\rr d\setminus 0)$ such that for every $\psi_1 \in \mathscr C^\dir_{x_0}(\rr d\setminus 0)$
and every $\psi_2 \in \mathscr C^\dir _{\xi _0}(\rr d\setminus 0)$, there exists $j\in J$ such that (and for every $j\in J$ it holds)
\begin{equation}\label{psiwfpsiecond3}
\psi _1\cdot \psi _2 (D)f\notin \cB_{j,k}
\end{equation}
for every $k\in J$ (for some $k\in J$).
\end{enumerate}

\noindent
Finally, the \emph{global wave-front sets of $\cupap$- and $\capup$-types} $\WF{\cupap}_{(\mathcal B
_{j,k})}(f)\subseteq (\rr d\times \rr d)\back 0$ and $\WF{\capup }_{(\mathcal
B _{j,k})}(f)\subseteq (\rr d\times \rr d)\back 0$ are the sets
\begin{equation*}
\WF{\cupap }_{(\mathcal B _{j,k})}(f) \equiv \WF {\psi ,\, \cupap}_{(\mathcal B
_{j,k})}(f)\bigcup \WF {e ,\, \cupap}_{(\mathcal B _{j,k})}(f) \bigcup \WF {\psi e ,\,
\cupap}_{(\mathcal B _{j,k})}(f),
\end{equation*}
and
\begin{equation*}
\WF{\capup}_{(\mathcal B _{j,k})}(f) \equiv \WF {\psi ,\, \capup}_{(\mathcal B
_{j,k})}(f)\bigcup \WF {e ,\, \capup}_{(\mathcal B _{j,k})}(f) \bigcup  \WF {\psi e
,\, \capup}_{(\mathcal B _{j,k})}(f)
\end{equation*}
respectively.
\end{defn}

\par

\begin{rem}
In analogy with Remark \ref{supinfequalwf} we note that if $\cB_{j,k}=\cB _j$ is independent of  $k\in J$, then
$$
\WF {\psi ,\, \cupap}_{(\mathcal B _{j,k})}(f)=\WF {\psi ,\, \cup}_{(\mathcal B _{j})}(f)
$$
and
$$
\WF {\psi ,\, \capup}_{(\mathcal B _{j,k})}(f)=\WF {\psi ,\, \cap}_{(\mathcal B _{j})}(f),
$$
and similarly for $\WF {e ,\, \cupap}_{(\mathcal B _{j,k})}(f)$, $\WF
{\psi e,\, \cupap}_{(\mathcal B _{j,k})}(f)$, $\WF {e ,\, \capup}_{(\mathcal B _{j,k})}(f)$ and $\WF
{\psi e,\, \capup}_{(\mathcal B _{j,k})}(f)$.
Hence, the family of wave-front sets in Definition \ref{defsuperposWF2} contain the wave-front sets in Definition \ref{defsuperposWF}.
\end{rem}

\par

\begin{rem}
We note that if $\mathcal B_{j,k}$ is $\SG$-admissible for every $j, k$ and $f\in \mathscr S$, then 
$$\WF {\psi ,\, \cupap}_{(\mathcal B _{j,k})}(f), \quad \WF {e ,\, \cupap}_{(\mathcal B _{j,k})}(f) \quad \text{and} \quad \WF {\psi e ,\, \cupap}_{(\mathcal B _{j,k})}(f)
$$ are closed subsets of $\rr d\times (\rr d \setminus 0)$, $(\rr d \setminus 0)\times \rr d$ and $(\rr d \setminus 0) \times (\rr d\setminus 0)$ respectively.
\end{rem}

\par

The following two results follows by similar arguments as in the proof of Proposition \ref{capcupprop1} and Theorem \ref{mainthm4}$'$. The details are left for the reader. 

\par

\renewcommand{\rubrik}{Proposition \ref{capcupprop1}$'$}

\begin{tom}

Let $\cB _{j,k}$ be the same as in Definition \ref{defsuperposWF2}, and let $f\in \mathscr S'(\rr d)$. Then
\begin{equation*}
\begin{alignedat}{2}
\bigcup _{j} \Big (\bigcap _{k}\WF {\psi }_{\mathcal B _{j,k}}(f) \Big ) &\subseteq \WF {\psi ,\, \cupap}_{(\mathcal B _{j,k})}(f),& \qquad
\bigcup _{j} \Big (\bigcap _{k} \WF {e }_{\mathcal B _{j,k}}(f) \Big ) &\subseteq \WF {e ,\, \cupap}_{(\mathcal B _{j,k})}(f),
\\[1ex]
\bigcup _{j} \Big (\bigcap _{k} \WF {\psi e}_{\mathcal B _{j,k}}(f) \Big ) &\subseteq \WF {\psi e,\, \cupap}_{(\mathcal B _{j,k})}(f),&
\qquad
\bigcup _{j} \Big ( \bigcap _{k} \WF {}_{\mathcal B _{j,k}}(f) \Big ) &\subseteq \WF {\cupap}_{(\mathcal B _{j,k})}(f),
\end{alignedat}
\end{equation*}
and
\begin{equation*}
\begin{alignedat}{2}
\bigcap _{j} \Big (\bigcup _{k}\WF {\psi }_{\mathcal B _{j,k}}(f) \Big ) &\subseteq \WF {\psi ,\, \capup}_{(\mathcal B _{j,k})}(f),& \qquad
\bigcap _{j} \Big ( \bigcup _{k} \WF {e }_{\mathcal B _{j,k}}(f) \Big ) &\subseteq \WF {e ,\, \capup}_{(\mathcal B _{j,k})}(f) 
\\[1ex]
\bigcap _{j} \Big ( \bigcup _{k} \WF {\psi e}_{\mathcal B _{j,k}}(f) \Big ) &\subseteq \WF {\psi e,\, \capup}_{(\mathcal B _{j,k})}(f),&
\qquad
\bigcap _{j} \Big ( \bigcup _{k} \WF {}_{\mathcal B _{j,k}}(f) \Big ) &\subseteq \WF {\capup}_{(\mathcal B _{j,k})}(f),
\end{alignedat}
\end{equation*}
\end{tom}

\par

\renewcommand{\rubrik}{Theorem \ref{mainthm4}$''$}

\begin{tom}
Let $\cB _{j,k}$ be $\SG$-admissible for every $j$ and $k$, and let $f\in\cS^\prime (\rr d)$. Then
\begin{alignat*}{5}
f&\in \bigcup _j\Big ( \bigcap _k \cB _{j,k} \Big ) &\quad &\Longleftrightarrow &\quad \bigcap _j\Big (\bigcup _k \WF {}_{\cB _{j,k}}(f) \Big ) &= \emptyset &\quad &\Longleftrightarrow &\quad \WF {\capup}_{(\cB _{j,k})}(f) &= \emptyset ,
\intertext{and}
f&\in \bigcap _j\Big (\bigcup _k \cB _{j,k} \Big ) &\quad &\Longleftrightarrow &\quad \bigcup _j \Big ( \bigcap _k \WF {}_{\cB _{j,k}}(f) \Big ) &= \emptyset &\quad &\Longleftrightarrow &\quad \WF {\cupap}_{(\cB _{j,k})}(f) &= \emptyset .
\end{alignat*}
\end{tom}

\par

\begin{rem}\label{remhormwfsets2}
Let $Q_0(\rr d)$ be the set of all $f\in C^\infty (\rr d)$ such that for some $N\ge 0$ (depending on $f$) we have $f^{(\alpha)}/\sigma _N\in L^\infty$, for every multi-index $\alpha$. We also let $Q(\rr d)$ be the set of all $f\in C^\infty (\rr d)$ such that for every multi-index $\alpha$, there exists an $N\ge 0$ (depending on $f$ and $\alpha$) such that $f^{(\alpha)}/\sigma _N\in L^\infty$.
In a similar way as in Remark \ref{remhormwfsets}, we may construct wave-front sets with respect to these spaces.

\par

In fact, let $p_{j,k},q_{j,k} \in [1,\infty]$, $\mathscr B_{j,k}=L_1^{p_{j,k},q_{j,k}}(\rr {2d})$,
$\omega _{j,k}(x,\xi ) = \eabs x^{-j}\eabs \xi^{k}$,  $\mathcal B_{j,k}$ be as
in \eqref{Jsetdef} for $j,k\in J=\mathbf N_0$ and set $\cC _{j,k}=\cB _{k,j}$. By similar arguments in \cite[Remark 2.18]{HTW} it follows that
$$
Q_0(\rr d)=\bigcup _j\Big ( \bigcap _k \cB _{j,k} \Big ),\quad Q(\rr d)=\bigcap _j\Big ( \bigcup _k  \cC _{j,k}\Big ) .
$$

\par

Now we let the wave-front sets of $\psi$-type with respect to $Q_0(\rr d)$ and $Q(\rr d)$ be defined by the formulas
$$
\WF {\psi}_{Q_0} (f) = \WF {\psi ,\capup }_{(\cB _{j,k})}(f),\quad \WF {\psi}_{Q} (f) = \WF {\psi ,\cupap }_{( \cC _{j,k})}(f).
$$
In the same way the wave-front sets of $e$-type, $\psi e$-type and the global wave-front sets with respect to $Q_0(\rr d)$ and $Q(\rr d)$ are defined.

\par

If $f\in \mathscr S'(\rr d)$, then by Theorem \ref{mainthm4}$''$ it follows that
\begin{alignat*}{3}
f &\in Q_0(\rr d)& \qquad &\Longleftrightarrow & \qquad \WF {}_{Q_0}(f) &= \emptyset
\intertext{and}
f &\in Q(\rr d)& \qquad &\Longleftrightarrow & \qquad \WF {}_{Q}(f) &= \emptyset .
\end{alignat*}
\end{rem}

\par

\subsection{Mapping properties for pseudo-differential operators}

\par

Next we consider mapping properties for pseudo-differential operators on wave-front sets of sequence types. The following result follows immediately from the definitions, Theorem \ref{mainthm2} and its proofs.

\par

\renewcommand{\rubrik}{Theorem \ref{mainthm2}$'$}

\begin{tom}
Let $r,\rho \in [0,1]$, $t\in \mathbf R$, $\omega _0\in \mathscr P_{r,\rho }
(\rr {2d})$, $a\in \SG^{(\omega _0)}_{r,\rho} (\rr {2d})$  and let $f\in
\mathscr S'(\rr d)$. Moreover, let $(\cB _{j,k}, \cC _{j,k})$ be a $\SG$-ordered pair with respect to $\omega _0$ for every $j,k\in J$.
Then the following is true:
\begin{enumerate}
\item  if in addition $\rho >0$, then
\begin{multline}\tag*{(\ref{pwavefrontemb1})$'$}
    \WF{\psi ,\, \cupap}_{(\mathcal C_{j,k})} ( \op _t(a)f) \subseteq
    \WF{\psi ,\, \cupap}_{(\mathcal B_{j,k})} (f)
\\[1ex]
          \subseteq
          \WF{\psi ,\, \cupap}_{(\mathcal C_{j,k})} ( \op _t(a)f) \bigcup
\Char ^{\psi }_{(\omega _0)}(a)\text ;
\end{multline}

\vrum

\item  if in addition $r >0$, then
\begin{multline}\tag*{(\ref{ewavefrontemb1})$'$}
    \WF{ e,\, \cupap}_{(\mathcal C_{j,k})} ( \op _t(a)f) \subseteq
    \WF{ e,\, \cupap}_{(\mathcal B_{j,k})} (f)
\\[1ex]
          \subseteq
          \WF{e ,\, \cupap}_{(\mathcal C_{j,k})} ( \op _t(a)f)\bigcup
\Char ^{e}_{(\omega _0)}(a)\text ;
\end{multline}

\vrum

\item  if in addition $r,\rho >0$, then
\begin{multline}\tag*{(\ref{pewavefrontemb1})$'$}
    \WF{\psi e,\, \cupap}_{(\mathcal C_{j,k})} ( \op _t(a)f) \subseteq
    \WF{\psi e,\, \cupap}_{(\mathcal B_{j,k})} (f)
\\[1ex]
          \subseteq
          \WF{\psi e,\, \cupap}_{(\mathcal C_{j,k})} ( \op _t(a)f) \bigcup
\Char ^{\psi e}_{(\omega _0)}(a) .
\end{multline}
\end{enumerate}
The same is true if the wave-front sets of $\cupap$-types are replaced by
$\capup$-types.
\end{tom}

\par

We note that several properties that are valid for the wave-front
sets of modulation space types also hold for wave-front sets in the
present section. 
The following generalization of Theorem \ref{mainthm4} is an immediate
consequence of Theorem \ref{hypoellthm}$'$, since $\Char_{(\omega
_0)}(a) =\emptyset$, when $a$ is $\SG$-elliptic with respect to
$\omega_0$.

\par

\renewcommand{\rubrik}{Theorem \ref{hypoellthm}$'$}

\begin{tom}
Let $r,\rho \in (0,1]$, $t\in \mathbf R$, $\omega _0\in \mathscr P_{r,\rho}(\rr {2d})$ and let $a\in
\SG^{(\omega _0)}_{r,\rho}(\rr {2d})$ be $\SG$-elliptic with respect to $\omega
_0$ and let $f\in \cS '(\rr d)$. Moreover, let $(\cB _{j,k}, \cC _{j,k})$ be a $\SG$-ordered pair with respect to $\omega _0$ for every $j,k\in J$. Then
\begin{align*}
\WF{\psi ,\, \cupap} _{(\mathcal C _{j,k})} (\op _t(a)f) &= \WF{\psi ,\,
\cupap}_{(\mathcal B _{j,k})} (f),
\\[1ex]
\WF{e,\, \cupap}_{(\mathcal C _{j,k})}
(\op _t(a)f) &= \WF{e,\, \cupap}_{(\mathcal B _{j,k})} (f)
\\[1ex]
\WF{\psi e,\, \cupap}_{(\mathcal C _{j,k})} (\op _t(a)f) &= \WF{\psi e,\,
\cupap}_{(\mathcal B _{j,k})} (f).
\end{align*}
The same is true if the wave-front sets of $\cupap$-types are replaced by
$\capup$-types.
\end{tom}

\medspace

Next we list some consequences of the previous results. We are especially focused on mapping and wave-front properties in the framework of the spaces $\mathscr S$, $Q_0$ and $Q$.

\par

We start with the following result, where the first part is a slight extension of \cite[Theorem 1.1]{CoMa}. 

\par

\begin{prop}\label{thmclassicwfs1}
Let $r,\rho \in (0,1]$, $t\in \mathbf R$, and let $\omega _0\in \mathscr P_{\rho
,\delta }(\rr {2d})$, $a\in \SG ^{(\omega _0)}_{r,\rho}(\rr {2d})$,
$f\in \mathscr S'(\rr d)$. Then
\begin{alignat*}{2}
\WF{\psi} ( \op _t(a)f) &\subseteq  & \WF{\psi} (f) &\subseteq
          \WF{\psi}( \op _t(a)f) \bigcup \Char ^{\psi}_{(\omega _0)}(a),
\\[1ex]
\WF{e} ( \op _t(a)f) &\subseteq  & \WF{e} (f) &\subseteq
          \WF{e}( \op _t(a)f) \bigcup \Char ^{e}_{(\omega _0)}(a),
\\[1ex]
\WF{\psi e} ( \op _t(a)f) &\subseteq  & \WF{\psi e} (f) &\subseteq
          \WF{\psi e}( \op _t(a)f) \bigcup \Char ^{\psi e}_{(\omega _0)}(a).
\end{alignat*}

\par

The same is true if the wave-front sets with respect to $\mathscr S$ are replaced by wave-front sets with respect to $Q_0$ or $Q$.
\end{prop}

\par

\begin{proof}
The result is an immediate consequence of Remark \ref{remhormwfsets}, Theorems
\ref{mainthm4}$'$ and  \ref{hypoellthm}$'$ and Corollary
\ref{corhypoellthm}.
\end{proof}

\par

A combination of the previous result and Theorem \ref{hypoellthm}$'$ gives the following.

\par

\begin{prop}\label{thmclassicwfs2}
Let $0<r,\rho \le 1$, $t\in \mathbf R$, $\omega _0\in \mathscr P_{r,\rho}(\rr {2d})$ and
let $a\in \SG^{(\omega _0)}_{r,\rho}(\rr {2d})$ be $\SG$-elliptic with
respect to $\omega _0$. If $f\in \mathscr S'(\rr d)$, then
\begin{gather*}
\WF{\psi}  (\op _t(a)f) = \WF{\psi} (f),\quad \WF{e} (\op _t(a)f) =
\WF{e} (f),
\\[1ex]
\WF{\psi e} (\op _t(a)f) = \WF{\psi e} (f).
\end{gather*}

\par

The same is true if the wave-front sets with respect to $\mathscr S$ are replaced by wave-front sets with respect to $Q_0$ or $Q$.
\end{prop}

\par

The next result follows immediately from the previous proposition and Remark \ref{remhormwfsets2}.

\par

\begin{prop}\label{corhypoellthm2}
Assume that the hypothesis in Theorem \ref{hypoellthm}${}^{\, \prime}$
is fulfilled, and let $f,g\in \mathscr S(\rr d)$ be such that the equation
\begin{equation}\label{pseudoeq}
\op _t(a)f = g
\end{equation}
holds. Then the following is true:
\begin{enumerate}
\item $f\in \mathscr S$, if and only if $g\in \mathscr S$;

\vrum

\item $f\in Q_0$, if and only if $g\in Q_0$;

\vrum

\item $f\in Q$, if and only if $g\in Q$.
\end{enumerate}
\end{prop}

\par

\begin{example}
Let $r,\rho$, $a$ and $\omega_0$ be the same as in Example \ref{examplheat}.
%
%
%
%
Then we may apply Theorems \ref{hypoellthm}$'$ and Propositions \ref{thmclassicwfs2} and \ref{corhypoellthm2} on $\op (a)$. In particular, (1)--(3) in Proposition \ref{corhypoellthm2} holds when \eqref{pseudoeq} is fulfilled.
%
%
\end{example}

\begin{example} Let $x=(x _1, x _2)\in \rr d$ and $\xi=(\xi _1, \xi _2)\in\rr d$, where $x _1, \xi _1 \in V _1$ and $x _2, \xi _2\in V _2$  and $V _1 \oplus V _2=\rr d$. We may consider wave-front sets with respect a space of distributions which behave like $\mathscr S$ in the $x _1$-variable and like $\mathscr S'$ in the $x _2$-variable. 
More precisely, let
\begin{align*}
\omega _{j,0}(x_1,\xi _1) = \eabs {x _1,\xi _1}^{j},\quad \omega _{0,k}(x_2,\xi _2) &= \eabs {x _2,\xi _2}^{-k}
\\[1ex]
\text{and}\quad \omega_{j,k}(x,\xi ) = \scal {x _1}{\xi _1}^j\scal {x _2}{&\xi _2}^{-k}.
\end{align*}
Then
$$
\bigcap_{j\geq 0}M^{p,q}_{(\omega_{j,0})}(V_1)=\mathscr S(V _1) \qquad \text{and} \qquad \bigcup_{k\geq 0} M^{p,q}_{(\omega_{0,k})}(V_2)=\mathscr S'(V _2) .
$$
Hence, we may interprete the set
$$
\mathcal B=\bigcap_{j\geq 0}\big(\bigcup_{k\geq 0} M^{p,q}_{(\omega_{j,k})}(\rr d)\big)
$$
as $\mathscr S(V_1;\mathscr S'(V_2))$, the set of all tempered distributions which behaves like $\mathscr S$ in the $x_1$-variable, and like $\mathscr S'$ in the $x _2$-variable.

\par

The wave-front set with respect to $\mathcal B$ is the wave-front set of $\cupap$-type with respect to the sequence of modulation spaces $(M^{p,q}_{(\omega_{j,k})}(\rr d))_{j,k\geq 0}$. In particular, all the mapping properties (e.{\,}g. Theorems \ref{mainthm2}$'$ and \ref{hypoellthm}$'$) hold for such wave-front sets.
\end{example}

\par

The set of characteristic points which we considered so far are defined in terms of certain elliptic conditions for the given symbol class. In what follows we give examples on how the results in the present and previous section can be extended by replacing these elliptic types of set of characteristic points with hypoelliptic ones. For this reason, let $r,\rho \in [0,1]$, and let $\omega _1,\omega _2\in \mathscr P_{r,\rho}(\rr {2d})$ be such that
\begin{equation}\label{hypoellcond1}
C^{-1}\eabs x^{-\delta r}\eabs \xi ^{-\delta \rho}\omega _2(x,\xi )\le \omega _1(x,\xi ) \le C\omega _2(x,\xi ) ,
\end{equation}
for some constants $C>0$ and $\delta <1$. Then $a\in \SG ^{(\omega _2)}_{r,\rho}(\rr {2d})$ is called $\SG$-hypoelliptic with respect to $(\omega _2,\omega _1)$, if there is a constant $c>0$ such that
\begin{equation}\label{hypoellcond2}
c\omega _1 (x,\xi )\le |a(x,\xi )| 
\end{equation}
outside a compact set in $\rr {2d}$.

\par

The symbol $a$ is called $\psi$-hypoelliptic at $(x_0,\xi _0)\in \rr d\times (\rr d\setminus 0)$, if for some neighbourhood $X$ of $x_0$, some conical neighbourhood $\Gamma$ of $\xi _0$ and some constant $R>0$, we have that \eqref{hypoellcond2} holds when $x\in X$, $\xi \in \Gamma$ and $|\xi |>R$. The set of characteristic points of $\psi$-type for $a$ with respect to $(\omega _2,\omega _1)$ is denoted by
\begin{equation}\label{extcharpsiset}
\Char ^\psi _{(\omega _2,\omega _1)}(a),
\end{equation}
and consists of all $(x_0,\xi _0)\in \rr d\times (\rr d\setminus 0)$, where $a$ fails to be $\psi$-hypoelliptic. The sets of characteristic points
\begin{equation}\label{extcharpsieset}
\Char ^e _{(\omega _2,\omega _1)}(a)\subseteq (\rr d\setminus 0)\times \rr d\quad \text{and}\quad \Char ^{\psi e} _{(\omega _2,\omega _1)}(a)\subseteq (\rr d\setminus 0)\times (\rr d\setminus 0)
\end{equation}
are defined in analogous ways. We observe that $\SG$-hypoellipticity of $a$ is equivalent to the fact that
$$
\Char ^\psi _{(\omega _2,\omega _1)}(a),\quad \Char ^e _{(\omega _2,\omega _1)}(a)\quad \text{and}\quad \Char ^{\psi e} _{(\omega _2,\omega _1)}(a)
$$
are empty.

\par

Now assume that $\omega _1,\omega _2\in \mathscr P_{r,\rho}(\rr {2d})$ fulfill \eqref{hypoellcond1} for some $\delta <1$. Then the same arguments as in the proof of Proposition \ref{psiecharequiv} in Appendix A, show that Propositions \ref{psicharequiv}, \ref{echarequiv} and \ref{psiecharequiv} hold for $\omega _0=\omega _2$, after each assumption $b\in \SG ^{(1/\omega _0)}_{r,\rho}(\rr {2d})$ has been replaced by $b\in \SG ^{(1/\omega _1)}_{r,\rho}(\rr {2d})$, and the sets of characteristic points have been replaced by the corresponding ones in \eqref{extcharpsiset} and \eqref{extcharpsieset}.

\par

The following extension of Theorem \ref{mainthm2}$'$ now follows by similar arguments as in the proof of Theorem \ref{mainthm2}. The details are left for the reader.

\par

\renewcommand{\rubrik}{Theorem \ref{mainthm2}$''$}

\begin{tom}
Let $r,\rho \in [0,1]$, $t\in \mathbf R$, $\omega _1,\omega _2\in \mathscr P_{r,\rho }
(\rr {2d})$ be such that \eqref{hypoellcond1} holds for some constants $C>0$ and $\delta <1$, $a\in \SG^{(\omega _2)}_{r,\rho} (\rr {2d})$  and let $f\in
\mathscr S'(\rr d)$. Moreover, let $(\cB _{j,k}, \cC _{j,k})$ and $(\cB _{j,k}, \ccD _{j,k})$ be a $\SG$-ordered pairs with respect to $\omega _2$ and $\omega _1$ respectively for every $j,k\in J$.
Then the following is true:
\begin{enumerate}
\item  if in addition $\rho >0$, then
\begin{multline}\tag*{(\ref{pwavefrontemb1})$''$}
    \WF{\psi ,\, \cupap}_{(\mathcal C_{j,k})} ( \op _t(a)f) \subseteq
    \WF{\psi ,\, \cupap}_{(\mathcal B_{j,k})} (f)
\\[1ex]
          \subseteq
          \WF{\psi ,\, \cupap}_{(\mathcal D_{j,k})} ( \op _t(a)f) \bigcup
\Char ^{\psi }_{(\omega _2,\omega _1)}(a)\text ;
\end{multline}

\vrum

\item  if in addition $r >0$, then
\begin{multline}\tag*{(\ref{ewavefrontemb1})$''$}
    \WF{ e,\, \cupap}_{(\mathcal C_{j,k})} ( \op _t(a)f) \subseteq
    \WF{ e,\, \cupap}_{(\mathcal B_{j,k})} (f)
\\[1ex]
          \subseteq
          \WF{e ,\, \cupap}_{(\mathcal D_{j,k})} ( \op _t(a)f)\bigcup
\Char ^{e}_{(\omega _2,\omega _1)}(a)\text ;
\end{multline}

\vrum

\item  if in addition $r,\rho >0$, then
\begin{multline}\tag*{(\ref{pewavefrontemb1})$''$}
    \WF{\psi e,\, \cupap}_{(\mathcal C_{j,k})} ( \op _t(a)f) \subseteq
    \WF{\psi e,\, \cupap}_{(\mathcal B_{j,k})} (f)
\\[1ex]
          \subseteq
          \WF{\psi e,\, \cupap}_{(\mathcal D_{j,k})} ( \op _t(a)f) \bigcup
\Char ^{\psi e}_{(\omega _2,\omega _1)}(a) .
\end{multline}
\end{enumerate}
The same is true if the wave-front sets of $\cupap$-types are replaced by
$\capup$-types.
\end{tom}

\par

\section*{Appendix A}

In this appendix we give a proof of Proposition
\ref{psiecharequiv}.

\par

\begin{proof}[Proof of Proposition \ref{psiecharequiv}]
We may assume that $t=0$ in view of Proposition \ref{psicontprop}. The equivalence between (1) and (2) follows easily by letting
$b(x,\xi )=\psi _1(x)\psi _2(\xi )/a(x,\xi )$, for appropriate $\psi _1\in \mathscr C_{x_0}^{\dir}$ and
$\psi _2\in \mathscr C_{\xi _0}^{\dir}$.

\par

(4) $\Rightarrow$ (3) is obvious in view of Remark \ref{psiinvremark}. Assume that (3) holds. We claim that
\begin{align}
|a(x,\xi )b(x,\xi )| &\ge 1/2
\tag*{(A.1)}
\intertext{holds when}
(x,\xi ) &\in \Gamma _1\times \Gamma _2,\ |x|\ge R,\ |\xi |\ge
R,
\tag*{(A.2)}
\end{align}
for some choice of conical neighbourhoods $\Gamma _1$
and $\Gamma _2$ of $x_0$ and $\xi _0$, respectively, and some $R>0$.
In fact, the assumptions imply that $ab=c+h_1$ for some
$h_1\in \SG ^{-r,0}_{r,\rho}+\SG ^{0,-\rho}_{r,\rho}$. By choosing $R$
large enough and $\Gamma _1$ and $\Gamma _2$ sufficiently small
conical neighbourhoods of $x_0$ and $\xi _0$, respectively, it follows
that $c(x,\xi )=1$ and $|h(x,\xi )|\le 1/2$ when (A.2)
holds. This gives (A.1). Since $|b|\le C/\omega$, it
follows that \eqref{invcond} is fulfilled, and (1) follows.

\par

It remains to prove that (2) implies (4). Assume therefore that (2) is true. Let $\psi _{1,k}\in \mathscr
C^\dir  _{x_0}(\Gamma _1)$ and $\psi _{2,k}\in C^\dir _{\xi _0}(\Gamma _2)$ for $k=1,\dots ,4$, be chosen such that
$$
b_1(x,\xi )\equiv \psi _{1,1}(x) \psi _{2,1}(\xi ) /a(x,\xi )\in \SG^{(1/\omega _0)}_{r,\rho},
$$
and $\psi _{j,k}=1$ on $\supp \psi _{j,k+1}$.
If $c_1=\psi _{1,1}\otimes \psi _{2,1}\in \SG ^{0,0}_{r,\rho}$, then it follows that
\begin{equation}\tag*{(A.3)}
\op (b_j)\op (a)=\op (c_j)+\op (h_j)
\end{equation}
holds for $j=1$ and some $h_1\in \SG ^{-r,-\rho }_{r,\rho}$.

\par
For $j\ge 2$ we now define $\widetilde b_j\in \SG ^{(1/\omega
)}_{r,\rho}$ by the Neumann series
$$
\op (\widetilde b_j) = \sum _{k=0}^{j-1}(-1)^k\op (\widetilde r_k),
$$
where $\op (\widetilde r_k)=\op (h_1)^k\op (b_1)\in \op (\SG ^{(\sigma
_{-k\rho ,-k\rho}/\omega _0)}_{r,\rho})$. Then (A.3) gives
\begin{multline}\tag*{(A.4)}
\op (\widetilde b_j)\op (a) = \sum _{k=0}^{j-1}(-1)^k\op (h_1)^k\op
(b_1)\op (a)
\\[1ex]
=\sum _{k=0}^{j-1}(-1)^k\op (h_1)^k(\op (c_1)+\op (h_1))
\\[1ex]
=\op (c_1)+\op (\widetilde h_{1,j})+\op (\widetilde h_{2,j}),
\end{multline}
where
\begin{align}
\op (\widetilde h_{1,j}) &= (-1)^j\op (h_1)^{j}\in \op (\SG
^{-jr,-j\rho}_{r,\rho})\tag*{(A.5)}
\intertext{and}
\op
(\widetilde h_{2,j}) &= -\sum _{k=1}^{j-1}(-1)^k\op (h_1)^k\op
(1-c_1)\in \op (\SG ^{0,0}_{r,\rho}).\notag
\end{align}

\par

By asymptotic expansions it follows that
\begin{equation}\tag*{(A.6)}\label{h2tilde}
\op (\widetilde h_{2,j}) = -\sum _{k=1}^{j-1}(-1)^k\op (1-c_1)\op
(h_1)^k + \op (\widetilde h_{3,j}) + \op (\widetilde h_{4,j}),
\end{equation}
for some $\widetilde h_{3,j}\in \SG ^{-r,-\rho}_{r,\rho}$ which is
equal to zero on $\supp c_1$ and $\widetilde h_{4,j}\in \SG
^{-jr,-j\rho}_{r,\rho}$. Now let $c$, $b_j$ and $r_k$ be defined
by the formulae
\begin{gather*}
c(x,\xi ) = \psi _{1,3}(x)\psi _{2,3}(\xi ),\quad  \op (b_j)= \op
(c)\op (\widetilde b_j)\in \op (\SG ^{(1/\omega _0)}_{r,\rho}),
\\[1ex]
\op (r_k) = \op (c)\op (\widetilde r_k)\in \op (\SG ^{(\sigma
_{-k\rho ,-k\rho}/\omega _0)}_{r,\rho}).
\end{gather*}
Then
$$
\op (b_j) = \sum _{k=0}^{j-1}(-1)^k\op (r_k)
$$
and (A.4)--(A.6) give
\begin{multline*}
\op (b_j)\op (a) = \op (c)\op (c_1)+\op
(c)\op (\widetilde h_{1,j})
\\[1ex]
-\sum _{k=1}^{j-1}(-1)^k\op (c)\op (1-c_1)\op (h_1)^k + \op (c)\op
(\widetilde h_{3,j}) + \op (c)\op (\widetilde h_{4,j}).
\end{multline*}
Since $c_1=1$ and $\widetilde h_{3,j}=0$ on $\supp c$, and every element of
$\op(\SG ^{-\infty ,-\infty }_{r,\rho})$ maps continuously ${\mathscr S}^\prime$ to $\mathscr S$, we find
\begin{align*}
\op (c)\op (c_1) &= \op (c)\mod \op (\mathscr S),
\\[1ex]
\op (c)\op (\widetilde h_{1,j}) &\in \op (\SG ^{-jr,-j\rho}_{r,\rho}),
\\[1ex]
\sum _{k=1}^{j-1}(-1)^k\op (c)\op (1-c_1)\op (h_1)^k &\in  \op
(\mathscr S),
\\[1ex]
\op (c)\op (\widetilde h_{3,j}) &\in \op (\mathscr S),
\intertext{and}
\op (c)\op (\widetilde h_{4,j}) &\in \op (\SG ^{-jr,-j\rho}_{r,\rho}).
\end{align*}
Hence, (A.3) follows for $c_j=c$ and some $h_j\in \SG
^{-jr,-j\rho}_{r,\rho}$.
By choosing $b\in \SG ^{(1/\omega )}_{r,\rho}$ such that
$$
b\sim \sum r_k,
$$
the argument above shows that $\op (b)\op (a) =\op (c)+\op (h)$, with
$
h\in \SG ^{-\infty ,-\infty}_{r,\rho}=\mathscr S,
$
and (4) follows. The proof is complete.
\end{proof}

\par

\end{document}